\documentclass[preprint,12pt]{elsarticle}

\usepackage{amssymb}
\usepackage{amsmath}
\usepackage{subcaption}
\usepackage{xcolor}
\usepackage{array}

\journal{Communications in Computational Physics for peer review}
\date{}  
\begin{document}

\begin{frontmatter}




\title{Numerical Error Analysis of the Poisson Equation under RHS Inaccuracies in Particle-in-Cell Simulations}

\author[aff1,aff3]{Kai Zhang\corref{cor1}}
\ead{zhangkai@huas.edu.cn}
\cortext[cor1]{Corresponding author}

\author[aff1]{Xiao Tao}
\ead{TaoXiao3587@gmail.com}

\author[aff2]{Weizong Wang}
\ead{wangweizong@buaa.edu.cn}

\author[aff2]{Bijiao He}
\ead{hbj@buaa.edu.cn}

\affiliation[aff1]{organization={Department of Computer and Electric Engineering, Hunan University of Arts and Science},
	addressline={3150 Dongting Avenue}, 
	city={Changde},
	postcode={415000}, 
	state={Hunan},
	country={China}}

\affiliation[aff2]{organization={School of Aerospace Engineering, Beihang University},
	addressline={37 Xueyuan Road, Haidian District}, 
	postcode={100191}, 
	state={Beijing},
	country={China}}

\affiliation[aff3]{organization={Hunan Provincial Key Laboratory for Control Technology of Distributed Electric Propulsion Arcarft},
	city={Changde},
	postcode={415000}, 
	state={Hunan},
	country={China}}

\begin{abstract}
Particle-in-Cell (PIC) simulations rely on accurate solutions of the electrostatic Poisson equation, yet accuracy often deteriorates near irregular Dirichlet boundaries on Cartesian meshes. While much research has addressed discretization errors on the left-hand side (LHS) of the Poisson equation, the impact of right-hand-side (RHS) inaccuracies—arising from charge density sampling near boundaries in PIC methods—remains largely unexplored. This study analyzes the numerical errors induced by underestimated RHS values at near-boundary nodes when solving the Poisson equation using embedded boundary finite difference schemes with linear and quadratic treatments. Analytical derivations in one dimension and truncation error analyses in two dimensions reveal that such RHS inaccuracies modify local truncation behavior differently: they reduce the dominant truncation error in the linear scheme but introduce a zeroth-order term in the quadratic scheme, leading to larger global errors. Numerical experiments in one-, two-, and three-dimensional domains confirm these findings. Contrary to expectations, the linear scheme yields superior overall accuracy under typical PIC-induced RHS inaccuracies. A simple RHS calibration strategy is further proposed to restore the accuracy of the quadratic scheme. These results offer new insight into the interplay between boundary-induced RHS errors and discretization accuracy in Poisson-type problems.
\end{abstract}



\end{frontmatter}



\section{Introduction}
\label{sec:intro}
Particle-in-Cell\cite{birdsall2018plasma} (PIC) simulation has become an indispensable tool in studying plasma physics, where a core computational task is solving the electrostatic Poisson equation on a mesh. Most modern PIC codes\cite{messmer2004parallel,qiang2016efficient,khaziev2018hpic,zhang2020parallel,zhang2021plume} for simulating plasma flows adopt Cartesian meshes, owing to advantages such as straightforward implementation, efficient parallelization, and favorable computational performance. While Cartesian PIC codes can, in principle, accommodate arbitrarily shaped boundaries, solution accuracy often degrades in the presence of irregular Dirichlet boundaries. This degradation arises because mesh nodes rarely align perfectly with the boundary surfaces. For geometries with high symmetry---such as spherical or cylindrical shapes—these challenges can be alleviated by developing customized PIC codes using spherical\cite{hutchinson2002ion,hutchinson2003ion,hutchinson2004ion} or cylindrical meshes\cite{singh1994three,singh1997potential,singh2000enhanced}. However, such specialized codes are only applicable to a limited class of problems and typically require substantial development effort. Other approaches, such as finite-volume-based methods\cite{kuhn2021picfoam} or finite-element methods\cite{fivaz1998finite,han2021pife,bai2021implicit}, fundamentally alter the structure of standard PIC codes, introducing additional complexity and reduced computational efficiency.

There were many existing techniques as to improving the solution accuracy of the Poisson equation with irregular domain on Cartesian meshes, each providing unique way of reducing the truncation errors for the near-boundary nodes---defined as nodes less than one grid spacing away from boundary interfaces along any coordinate direction. For example, Mayo's method\cite{mayo1984fast,mayo1985fast,mayo1992fast} that employ the integral formations for the near-boundary nodes. While this method allows the use of fast Poisson solvers\cite{hockney1970potenitial}, the procedures require using numerical quadrature to evaluate the integrals followed by spline interpolations. In solving more general elliptic equations, the immersed interface method (IIM) proposed by LeVeque and Li\cite{leveque1994immersed,li2003overview} introduces an innovative six-point stencil for the near-boundary nodes, which necessitates determining appropriate coefficients for the stencil points by solving local linear systems. Both approaches have been proved to produce second-order accuracy\cite{beale2007accuracy}. The embedded boundary finite volume methods \cite{johansen1998cartesian,devendran2017fourth} are viable alternatives, with the scheme of Devendran \textit{et al.} \cite{devendran2017fourth} achieving fourth-order accuracy.

Among second-order accurate  methods, the embedded boundary finite difference method offers an attractive alternative\cite{jomaa2005embedded,jomaa2010shortley}. This method embeds irregular Dirichlet boundaries within a rectangular (or cuboidal in 3D) computational domain and modifies the Poisson equation's discretization near the boundary using either the linear\cite{collatz1933bemerkungen,gibou2002second} or the quadratic interpolation schemes\cite{chen1997simple,shortley1938numerical}. For interior grid nodes located at least one grid spacing away from the Dirichlet boundaries along any axial direction, the standard five-point stencil, known to yield second-order truncation error, is retained. The resulting discretization leads to a linear system to be solved numerically.

Extensive prior work has investigated the error behavior of the embedded boundary finite difference method. It is widely recognized that both the linear and quadratic schemes improve the solution accuracy. In particular, the linear scheme was proved to give second-order accuracy\cite{collatz1933bemerkungen,gibou2002second}, while the quadratic attains third-order accuracy at the boundaries and maintain second-order accuracy in the interior\cite{shortley1938numerical,wasow1955discrete,matsunaga2000superconvergence}. These two schemes have been successfully applied in various computational applications, including crystal growth\cite{li2011fast}, fluid dynamics\cite{sharma2015level}, tumor modeling\cite{hogea2006simulating}, etc. Recently, the use of such schemes has also seen increases in plasma simulations\cite{revel2017massive,liu2021parallel}. The linear scheme, in particular, had been considered favorable\cite{gibou2002second} for its simplicity and the symmetric linear systems it produces, which enables using efficient sparse matrix storage and matrix-vector algorithms, as well as iterative solvers such as the Preconditioned Conjugate Gradient (PCG) method. However, Jomaa and Macaskill's error analysis\cite{jomaa2005embedded} challenged the preference for the linear scheme. In one-dimensional problems they provided explicit numerical error expressions. Upon evaluation of these expressions, they  showed that the linear scheme, although exhibiting  second-order accuracy, has much larger error coefficients than the quadratic counterpart. For two-dimensional problems, they verified numerically that their 1-D error analysis holds for non-corner grid points. 

While the majority of existing studies focus on discretization errors arising from the numerical treatment of the left-hand side (LHS) of the Poisson equation, comparatively little attention has been paid to inaccuracies originating from the right-hand side (RHS). Such inaccuracies may stem from measurement errors, data processing procedures, or physical modeling assumptions, and can significantly affect the accuracy of the computed potential field.

Several studies in the literature have examined the Poisson equation with varied RHS values. Pan \textit{et al.} \cite{pan2016error,faiella2021error} investigated the pressure Poisson equation in the context of particle image velocimetry (PIV) for fluid dynamics. In their formulation, the RHS depends on experimentally measured velocity field data, which inevitably contains systematic errors across the domain. They derived error bounds for the reconstructed pressure field and numerically validated these bounds under the assumption of constant perturbations in the RHS at both interior and boundary points. However, their analysis did not account for irregular Dirichlet boundary effects, nor did it consider numerical errors arising from the solving the Poisson equation. Marques \textit{et al.} \cite{marques2011correction} developed the Correction Function Method (CFM) to handle interface jump conditions for constant-coefficient Poisson equations. As a finite difference approach, this method estimates ghost-node values that can be incorporated into the RHS, thereby enabling the use of standard solvers. Although the method modifies the RHS, its primary purpose is to improve the discretization of the LHS; hence, it does not address inaccuracies in the RHS itself. 

A critical gap in the literature concerns the effect of RHS inaccuracies in the electrostatic Poisson equation, particularly in the context of PIC simulations where such inaccuracies frequently occur at nodes adjacent to Dirichlet boundaries. The electrostatic Poisson equation has the form $\nabla^2\phi = -\rho/\epsilon_0$, where $\phi$ is the electric potential to be solved, $\rho$ is the charge density, and $\epsilon_0$ is the vacuum permittivity. In standard PIC procedures, $\rho$ is computed as the weighted average of nearby particle charges over a cell volume. While this sampling algorithm yields accurate evaluation of charge density at interior nodes, provided that the grid spacing is small and the population of simulated particles are large, it’s seriously flawed for nodes near Dirichlet boundaries, because the averaging process assumes a full cell volume even though the cell may get cut off by the boundary interfaces. As a result, near-boundary nodes’ $\rho$ values---and thus the corresponding RHS values---are systematically underestimated due to the Dirichlet boundary's obstruction, as charged particles do not exist beyond the boundaries. Consequently, numerical solutions of the Poisson equation from all the existing numerical methods suffer from various degree of accuracy degradations in the presence of such RHS inaccuracies. To our knowledge, this potentially significant source of error has not been systematically investigated. 

Despite the lack of systematic investigation into near-boundary RHS inaccuracies, a few approaches have been proposed to mitigate related effects. In two-dimensional cylindrical meshes, Cornet and Kwok\cite{cornet2007new} improved the particle weighting algorithm of PIC method in a multiple-grid system. More recently, Lv and Zhong\cite{lv2025paro} employed machine learning techniques to develop a surrogate Poisson solver that directly maps $\rho$ data to the potential $\phi$, aiming primarily to reduce statistical noise. However, these efforts do not directly address the aforementioned issue of near-boundary RHS inaccuracies coupled with irregular Dirichlet boundaries.

In this study, we analyze the numerical errors arising from solving the electrostatic Poisson equation with embedded irregular Dirichlet boundaries on a Cartesian mesh, using both the linear and quadratic boundary treatment schemes, under the typical condition of underestimated RHS values near boundaries as encountered in PIC simulations. Our objective is to reveal how these RHS inaccuracies influence both the local truncation errors and the global numerical solution errors for interior and near-boundary nodes.

To achieve this, we devised numerical procedures that emulate the RHS inaccuracies introduced by standard PIC charge density sampling algorithms near Dirichlet boundaries. We then performed error analyses analogous to those of Jomaa and Macaskill\cite{jomaa2005embedded} for both one- and two-dimensional cases, incorporating the effects of RHS inaccuracies. For 1-D problems, we derived explicit error expressions; for 2D problems, we focused on analyzing truncation error magnitudes. We also validated our analysis via Poisson problems with pre-defined exact solutions in 1D, 2D, and 3D domains. Finally, we evaluated the performance of a simple RHS calibration strategy that mitigates the RHS inaccuracies.

Surprisingly, both our analytical and numerical results demonstrated that the linear scheme outperforms the quadratic scheme in the presence of inaccurate RHS values. The linear scheme not only reduces errors at near-boundary nodes but also improves overall accuracy across the domain, contrary to expectations. We show that the underestimated RHS values effectively alter the local truncation errors at near-boundary nodes: they reduce the magnitude of the (zeroth-order) truncation error in the linear scheme while introducing a zeroth-order term into the (originally first-order) truncation error of the quadratic scheme. These changes in truncation errors at the near-boundary nodes directly translate to the numerical errors at the boundary, which propagate to the interior and significantly affect the global solution accuracy. Our findings suggest that, under realistic RHS inaccuracies typical of PIC simulations, the linear scheme is preferred over the quadratic scheme. However, using a simple RHS calibration strategy, the accuracy of the quadratic scheme and, similarly, other higher order methods could be recovered. Consequently, our findings extend beyond the PIC community, offering mathematical insights into the numerical error behavior of the Poisson equation.

The remainder of this paper is organized as follows. In Section 2, we introduce the numerical formulation of the Poisson equation and describe the embedded boundary discretization schemes, along with the modeling of RHS inaccuracies. Section 3 presents the error analysis for both 1D and 2D cases, including the derivation of explicit expressions in 1D and truncation error evaluations in 2D. In Section 4, we validate our analysis through numerical experiments in 1D, 2D, and 3D, and evaluate the effectiveness of a simple RHS calibration strategy. Finally, Section 5 summarizes the key findings and discusses implications for PIC simulations and broader applications.

\section{Mathematical Model and Numerical Formulation}
\label{sec:numerical}
\subsection{Poisson equation and discretization schemes}
\label{subsec:poisson}
The electrostatic Poisson equation is discretized into the standard 2D finite difference form
\begin{equation}
\label{eqn:standard}
    \frac{\phi_{i-1,j}-2\phi_{i,j}+\phi_{i+1,j}}{\Delta_x^2}+
    \frac{\phi_{i,j-1}-2\phi_{i,j}+\phi_{i,j+1}}{\Delta_y^2}=-\frac{\rho_{i,j}}{\epsilon_0}=b_{i,j},
\end{equation}
where $\phi_{i,j}$ denotes the electrostatic potential at grid node $(i,j)$, $\Delta_x$ and $\Delta_y$ are the grid spacing along the $x$ and $y$ directions, and $b_{i,j}$ is the computed RHS value. This formulation can be readily extended to 3D, or reduced to 1D, by adding or removing corresponding axial terms. Each interior node is associated with one such discrete equation. For nodes adjacent to Dirichlet boundaries, terms involving known boundary values $\phi^D$ are transferred to the RHS. Combining all these equations yields a set of linear equations of the form
\begin{equation}
\label{eqn:ax=b}
    \textbf{A}\Vec{\phi}=\Vec{b},
\end{equation}
where \textbf{A} is the coefficient matrix, $\Vec{\phi}$ is the vector of unknown potentials at interior nodes, and $\Vec{b}$ is the corresponding RHS values vector. The matrix A is symmetric positive definite (SPD), allowing the use of efficient iterative solvers. The resulting numerical solution is second-order accurate.

\begin{figure}
    \centering
    \includegraphics[width=0.65\linewidth]{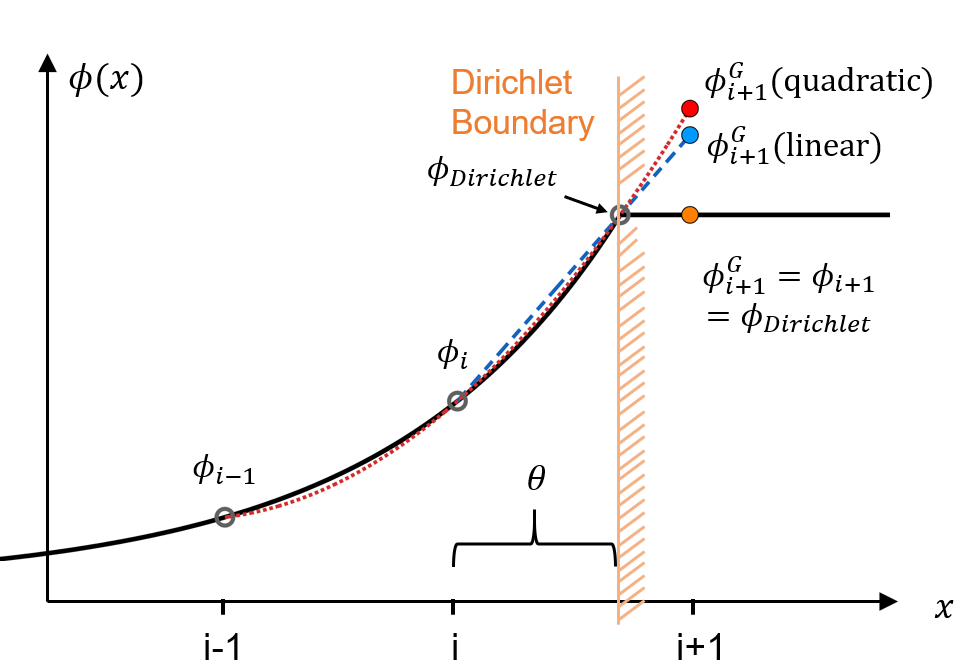}
    \caption{1D problem illustration with a Dirichlet boundary positioned between nodes.}
    \label{fig:1Dillustrate}
\end{figure}

The standard finite difference scheme performs well when the Dirichlet boundaries align with grid nodes; however, this condition is rarely met in practice, leading to a degradation of solution accuracy to first-order. Consider a 1D scenario where node $i+1$ lies outside the Dirichlet boundary. In this case, the discretization near node $i$ requires modification, as illustrated in figure~\ref{fig:1Dillustrate}. One approach is to use a linear extrapolation between node $(i)$ and the boundary point $\phi_D$ along the $x$-direction , resulting
\begin{equation}
\label{eqn:linear}
    \phi_{i+1}^G=\frac{\phi_D+(\theta-1)\phi_i}{\theta},
\end{equation}
where $\phi_{i,j}^G$ denotes a ghost value, and $\theta$ is the normalized distance from node $i$ to the boundary. An alternative is to apply a quadratic extrapolation using node $i,$ $i-1$, and the boundary point, resulting in
\begin{equation}
\label{eqn:quad}
    \phi_{i+1}^G=\frac{2}{\theta^2+\theta}\phi_D+\frac{2\theta-2}{\theta}\phi_i+\frac{1-\theta}{1+\theta}\phi_{i-1}.
    \end{equation}
Both modified schemes achieve second-order accuracy, though the linear version exhibits a larger error coefficient\cite{jomaa2005embedded}. A notable distinction is that the linear scheme preserves the symmetry of the coefficient matrix $\textbf{A}$, which can be advantageous for applying certain solvers.

\subsection{Numerical treatments for RHS values}
\label{subsec:rhs}
\begin{figure}
    \centering
    \includegraphics[width=0.7\linewidth]{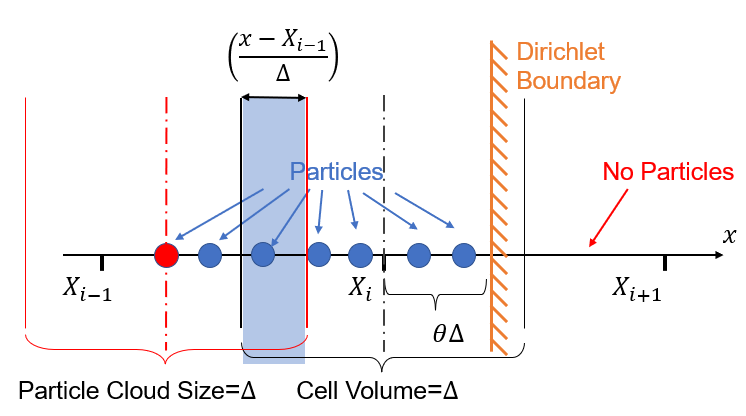}
    \caption{RHS inaccuracy encountered in PIC simulations, illustrated in 1D}
    \label{fig:rho_illust}
\end{figure}

As discussed earlier, the RHS values of the electrostatic Poisson equation in PIC simulations is directly determined by the charge density $\rho$, which is often computed inaccurately near irregular shaped Dirichlet boundaries. The charge density at an interior node is normally calculated as a weighted average of electric charges from surrounding particles, normalized by the cell volume. Figure~\ref{fig:rho_illust} illustrate this process in 1D setting. In PIC procedures, each particle represents a finite-sized "cloud" of charge with a volume equivalent to one cell, following the the cloud-in-cell model~\cite{birdsall2018plasma}. Charges from this cloud are distributed to nearby grid nodes using first-order weighting---equivalent to linear interpolation in 1D. For example, the red particle in figure~\ref{fig:rho_illust} deposits a fraction $(x-X_{i-1})/\Delta$ of its charge to node $i$, with the remainder assigned to node $i-1$. The total deposited charge at a node is then divided by the cell volume to yield the local charge density $\rho$. However, this approach becomes inaccurate near Dirichlet boundaries when a node lies within one grid spacing of the boundary. In such cases, part of the node’s effective deposition volume is truncated by the boundary interface—beyond which no particles exist—leading to an underestimation of both the charge density $\phi$ and the corresponding RHS value $b$. This local error can significantly affect the overall solution accuracy.

\begin{figure}
  \begin{subfigure}{0.49\textwidth}
    \includegraphics[width=\linewidth]{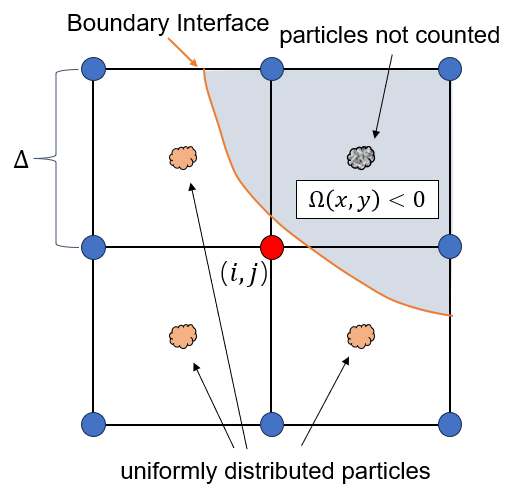}
    \caption{Particles distribution at level 1 partition} \label{fig:level1}
  \end{subfigure}%
  \hspace*{\fill}   
  \begin{subfigure}{0.49\textwidth}
    \includegraphics[width=\linewidth]{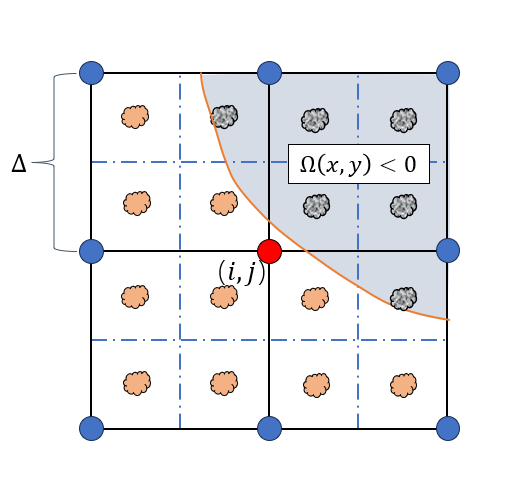}
    \caption{Particles distribution at level 2 partition} \label{fig:level2}
  \end{subfigure}

\caption{2D illustration for RHS value treatment for a near-boundary node $(i,j)$} \label{fig:1}
\label{fig:levels}
\end{figure}

To analyze the impact of these RHS inaccuracies on the overall solution quality, we designed a numerical procedure that emulates the typical RHS errors arising in PIC simulations at near-boundary nodes. For interior nodes free from RHS inaccuracies, their RHS values are computed directly from the exact solution using $\nabla^2\phi(\vec{x})$. For near-boundary nodes, however, we mimic the PIC routine by performing particle-based sampling of the surrounding region. The Dirichlet boundary interface is defined explicitly by a level-set function $\Omega(\vec{x})=0$, where $\Omega<0$ indicates the exterior and $\Omega\geq0$ defines the interior computational domain. Figure~\ref{fig:levels} illustrates this setup in 2D. For a near-boundary node $(i,j)$ (solid red circle), instead of directly applying the exact solution to compute its RHS, we populate nearby cells with virtual "particles" to replicate the charge deposition behavior. These particles are uniformly distributed at subcell centers, with the number determined by a partition level $l$. In 2D, each node is associated with $2^{2l}$ particles, and in 3D with $2^{3l}$ particles. For each particle $k$, the exact RHS value $b_k$ is evaluated at its location. Each particle contributes to node $(i,j)$'s RHS according to first-order weighting. Summing over all contributions yields the approximated RHS value $\overline{b_{i,j}}$, computed as:

\begin{equation}
    \label{eqn:rhs_calc}
    \overline{b_{i,j}}=\sum_{k=1}^{2^{Dl}}\frac{1}{2^{D(l-1)}}b_k\left( 1-\frac{|x_k-X_{i,j}|}{\Delta} \right) \left( 1-\frac{|y_k-Y_{i,j}|}{\Delta} \right),
\end{equation}
where the overline in $\overline{b_{i,j}}$ indicates that the value is approximated (inaccurate), $D$ is the spatial dimension (e.g., $D=2$ for 2D problems), $b_k$ is the exact RHS value at the $k$-th particle location, and $(X_{i,j},Y_{i,j})$ are the coordinates of node $(i,j)$. Only particles that lie within the computational domain ($\Omega$$\geq0$) are included in the summation, thereby replicating the physical absence of charge beyond the boundary. For 1D and 3D cases, the weighting terms in equation\eqref{eqn:rhs_calc} should be adjusted accordingly to accommodate the dimensionality. This treatment provides a controlled and physically consistent way to reproduce the RHS underestimation effect observed in standard PIC routines near boundaries.

 To facilitate subsequent error analysis, we define a parameter $\delta$ as the ratio between the inaccurate and accurate RHS values at a node, i.e., $\delta=\overline{b}/b$. Apparently, $\delta_{i,j}$ lies within the range $\in(0,1]$. In later numerical tests where the true RHS values $b_{i,j}$ are known, the corresponding $\delta_{i,j}$ values can be computed directly from this definition. However, in practical PIC simulations, the exact RHS distribution is not available a priori. Therefore, some assumption about the $b(\vec{x})$ distribution must be made to estimate an approximated $\overline{\delta_{i,j}}$. In our approach, we assume a locally uniform distribution of $b(\vec{x})$ in the vicinity of node $(i,j)$. The parameter $\delta_{i,j}$ (and its approximation $\overline{\delta}$) is implicitly tied to the boundary geometry, represented by the level-set function $\Omega(\vec{x})$, making direct computation difficult. Nevertheless, in the 1D scenario illustrated in figure~\ref{fig:rho_illust}, ${\delta_{i}}$ could be derived explicitly by the charge deposition integrals over the truncated cell volume. Specifically, we have
\begin{equation}
\label{eqn:1D_delta_int}
\delta_i=\overline{b_i}/{b_i}=\frac{\int_{X_{i-1}}^{X_i}{b(x)\cdot}\left(\frac{x-X_{i-1}}{\Delta}\right)dx+\int_{X_i}^{X_{i+\theta\Delta}}{b(x)\cdot}\left(\frac{X_{i+1}-x}{\Delta}\right)dx}{\Delta}/b_i,
\end{equation}
where $b(x)$ is an abstract (unknown) RHS distribution function, and the boundary truncates the cell at a distance $\theta\Delta$ from node $i$. Assuming $b(x)$ is locally uniform $b(x)$ yields a closed-form expression for the approximated ratio:
\begin{equation}
\label{eqn:1D_delta}
    \overline{\delta_i}=-\frac{1}{2}\theta^2+\theta+\frac{1}{2}.
\end{equation}
This expression is intuitive when considering the limiting cases. For $\theta=0$, the Dirichlet boundary lies exactly at node $i$, blocking half of the cell's volume. Accordingly, $\overline{\delta_i}=1/2$. For $\theta=1$, the boundary lies at the neighboring node $i+1$, and node $i$ behaves as a normal interior node with no volume loss, giving $\overline{\delta_i}=1$, as expected.

In the multidimensional (2D or 3D) setting, we extend the same idea used in the 1D case to compute the approximated coefficient $\overline{\delta_{i,j}}$ (or $\overline{\delta_{i,j,k}}$ in 3D) for near-boundary nodes. Specifically, we apply the same numerical procedure used in the RHS computation process (as described in equation\eqref{eqn:rhs_calc}), but under the assumption that $b(\vec{x})$ is locally uniform. For example, in 2D, we subdivide the four neighboring cells surrounding a near-boundary node $(i,j)$ into smaller subcells and place particles at the center of each subcell, following the same sampling strategy as in the actual RHS computation. Each particle is assigned a constant RHS value (i.e., $b_k = b$ for some arbitrary constant $b$), so the resulting computed value $\overline{b_{i,j}}$ represents only a fraction of the full uniform value. This fraction, by construction, corresponds to the approximated coefficient $\overline{\delta_{i,j}}$, and it can be evaluated using the expression
\begin{equation}
\label{eqn:delta}
    \overline{\delta_{i,j}}=\overline{b_{i,j}}/b=\sum_{k=1}^{2^{Dl}}\frac{1}{2^{D(l-1)}}\left( 1-\frac{|x_k-X_{i,j}|}{\Delta} \right) \left( 1-\frac{|y_k-Y_{i,j}|}{\Delta} \right),
\end{equation}
where $\overline{b_{i,j}}$ denotes the value computed using the assumed uniform distribution, and the summation spans all particle contributions within the neighboring region whose positions satisfy $\Omega \geq 0$, ensuring that only particles located inside the physical domain are counted. This construction allows us to numerically capture the influence of irregular geometry on the RHS evaluations at boundary-adjacent nodes, despite the lack of an explicit form for $\overline{\delta_{i,j}}$ in terms of $\theta_x$ and $\theta_y$.
\section{Error Analysis}
\subsection{1-D Error analysis}
We extend the earlier 1-D error analysis\cite{jomaa2005embedded} to account for inaccuracies of RHS values $\overline{b_i}$ in both the linear and quadratic embedded finite difference scheme. Consider 1-D Poisson equation posed on a domain $[a,b]$ with Dirichlet boundaries located at $x_L$ and $x_R$. The domain is discretized using $N+1$ uniformly spaced grid points such that $a=x_0<x_L<x_1<x_2<\cdot<x_{N-1}<x_R<b=x_N$. The fractional distances from $x_L$ and $x_R$ to their respective adjacent grid nodes are denoted by $\theta_L$ and $\theta_R$, with $x_L - x_1 = \theta_L \Delta_x$ and $x_{N-1} - x_{R} = \theta_R \Delta_x$.

Let the numerical error at node $i$ be defined as $\xi_i = \phi_i - \phi_i^e$, where $\phi_i$ is the numerical solution and $\phi_i^e$ is the exact solution. As shown by Jomaa and Macaskill\cite{jomaa2005embedded}, the numerical error of any grid node $\xi_i$ is an aggregate effect of all the truncation errors $\tau_i$, which is defined as:
\begin{equation}
	\tau_i=b_i-(L\phi^e)_i,
\end{equation}
where $L$ denotes the discrete Laplacian operator, which corresponds to the standard second-order central difference scheme for interior nodes. At the two near-boundary nodes $i = 1$ and $i = N-1$, $\overline{b}_i$ is used instead of the exact RHS, and $L$ must be adjusted based on whether the linear or quadratic embedded scheme is employed. A crucial relation between numerical error $\xi$ and truncation error $\tau$ could be established as: 
\begin{equation}
	\label{eqn:xiandtau}
	(L\xi)_i=\tau_i=\frac{H_{i+1/2}-H_{i-1/2}}{\Delta_x},\quad 1\leq i\leq N-1
\end{equation}
 where $H$ is a first-order difference operator acting on $\xi$, defined as
 \begin{equation}
 	H_{i-1/2}=\frac{\xi_i-\xi_{i-1}}{\Delta_x}.
 \end{equation}
 Note that at the domain boundaries, $H_{1/2}$ and $H_{N-1/2}$ are defined consistently with the chosen boundary discretization scheme used for $\phi$.
 
 \subsubsection{The linear scheme case}
 For interior nodes, the truncation error is obtained via Taylor expansion: 
 \begin{equation}
 	\tau_i=-\frac{\Delta_x^2}{12}\phi_i^{(4)}+\mathcal{O}(\Delta_x^4),\quad 2\leq i\leq N-2
 \end{equation}
which is a second-order term. Without considering RHS inaccuracies, the truncation errors at the near-boundary nodes are:
\begin{equation*}
	\tau_1=\frac{1}{2}(1-\theta_L)\phi_1^{''}+\mathcal{O}(\Delta_x),
\end{equation*} 
and
\begin{equation*}
	\tau_{N-1}=\frac{1}{2}(1-\theta_R)\phi_{N-1}^{''}+\mathcal{O}(\Delta_x),
\end{equation*}
both being zeroth-order.

When incorporating RHS inaccuracies and assuming a uniform $b(x)$ distribution, the truncation error at node $1$ becomes:
\begin{eqnarray*}
	\tau_1=\overline{b_1}-(L\phi^e)_1 &=& b_1-(L\phi^e)+\delta_1 b_1-b_1 \nonumber\\
	&=& \frac{1}{2}(1-\theta_L)\phi_1^{''}+(\delta_1-1)b_1+\mathcal{O}(\Delta_x).
\end{eqnarray*}
And by substituting $\delta_1$ from equation\eqref{eqn:1D_delta} and realizing $b_1=\phi_1^{''}$, we have
\begin{equation}
	\label{eqn:tau1_linear}
	\tau_1=\frac{\theta_L-\theta_L^2}{2}\phi_1^{''}+\mathcal{O}(\Delta_x),
\end{equation} 
which is still a zeroth-order term. Similarly, for the right boundary:
\begin{equation}
	\label{tau_n-1_linear}
	\tau_{N-1}=\frac{\theta_R-\theta_R^2}{2}\phi_{N-1}^{''}+\mathcal{O}(\Delta_x)
\end{equation}

The above relations enable us to write all the numerical error terms to $\xi_1$ as
\begin{equation}
	\label{eqn:xi_to_xi1}
	\xi_i=\xi_1+(i-1)\Delta_xH_{N-1/2}-\Delta_x^2\sum_{j=1}^{i-1}\sum_{k=j+1}^{N-1}\tau_k.\quad 2\leq i\leq N-1
\end{equation}
To close the set of equations, the left boundary treatment is needed, treated with the linear scheme by
\begin{equation}
	H_{1/2}=\frac{\xi_1-\xi_L}{\theta_L\Delta_x}=\frac{\xi_1}{\theta_L\Delta_x},
\end{equation}
where $\xi_L=\phi_L-\phi_L^e=0$ by definition. Also needed is the linear treatment at the right boundary
\begin{equation}
	H_{N-1/2}=\frac{\xi_R-\xi_{N-1}}{\theta_R\Delta_x}=\frac{-\xi_{N-1}}{\theta_R\Delta_x}.
\end{equation}
Solving these equations finally result in
\begin{equation}
	\label{eqn:xi_as_tau}
	\xi_i=\Delta_x^2\left[\left(\frac{i+\theta_L-1}{N+\theta_L+\theta_R-2}-1\right)\sum_{k=1}^{N-1}(k+\theta_L-1)\tau_k-\sum_{k=i+1}^{N-1}(i-k)\tau_k\right].
\end{equation} 
This explicit formula expresses the numerical error $\xi_i$ with all the truncation error terms.

To isolate the impact of of the left Dirichlet boundary, we set $\theta_R = 1$ and consider only $\tau_1$:
\begin{equation}
	\xi_i^L=\left(\frac{i}{N}-1\right)\theta_L\tau_1\Delta_x^2\simeq\left(\frac{i}{N}-1\right)\frac{\theta_L^2(1-\theta_L)}{2}\phi_1^{''}\Delta_x^2,
\end{equation}
where $\tau_1$ retained only the highest order term. Apparently, the left boundary induced error component is also a second-order in $\Delta_x$, linearly decreasing from left ro right. By taking the derivative we found that maximum error occurs $\theta_L=2/3$ (without RHS inaccuracies this value was found to be $1/2$\cite{jomaa2005embedded}). Similarly, for the right boundary we have
\begin{equation}
	\xi_i^R=-\frac{i}{N}\theta_R\tau_{N-1}\Delta_x^2\simeq-\frac{i}{N}\frac{\theta_R^2(1-\theta_R)}{2}\phi_{N-1}^{''}\Delta_x^2,
\end{equation}
which is also second-order in $\Delta_x$ and peaks at $\theta_R=2/3$. From the above expressions, we also observe that the boundary-induced error components $\xi_i^{L}$ and $\xi_i^R$ consistently have signs opposite to those of the corresponding boundary truncation errors $\tau_1$ and $\tau_{N-1}$.

The contribution from interior nodes can be approximated by the integral
\begin{equation}
	\label{eqn:inner_contrib}
	\begin{split}
		\xi_i^{in} &\simeq\frac{\Delta_x^2}{12} \Bigg[ \int_{x_L}^{x_i}x\phi^{(4)}(x)dx 
		-\frac{x_i-x_L}{x_R-x_L}\int_{x_L}^{x_R}x\phi^{(4)}(x)dx \\
		&\quad + x_L\frac{x_i-x_L}{x_R-x_L}\int_{x_L}^{x_R}\phi^{(4)}(x)dx-x_L\int_{x_L}^{x_i}\phi^{(4)}(x)dx\\
		&\quad + (x_i-x_L)\int_{x_i}^{x_R}\phi^{(4)}(x)dx 
		\Bigg],
	\end{split}
\end{equation}
which is also second-order in $\Delta_x$. Thus, the total numerical error at node $i$ can be approximated as: $\xi_i \approx \xi_i^L + \xi_i^R + \xi_i^{in}$. This straightforward error decomposition offers a clear and intuitive interpretation of the distinct error contributions in the numerical solution of the Poisson equation.

\subsubsection{The quadratic scheme case}
Adopting a quadratic scheme improves the truncation errors at the boundaries only, while the truncation errors at interior nodes remain unchanged. In the absence of the inaccurate RHS effect, the truncation errors at the near-boundary nodes $\tau_1$ and $\tau_{N-1}$ are given by
\begin{equation*}
	\tau_1=-\frac{(1-\theta_L)}{3}\phi_1^{'''}\Delta_x+\mathcal{O}(\Delta_x^2),
\end{equation*} 
and
\begin{equation*}
	\tau_{N-1}=\frac{(1-\theta_R)}{3}\phi_{N-1}^{'''}\Delta_x+\mathcal{O}(\Delta_x^2),
\end{equation*}
both of which are first-order in $\Delta_x$.

Similar to the analysis for the linear scheme, when the RHS inaccuracies are introduced under the assumption of a locally uniform RHS, the boundary truncation errors are significantly altered. For the left boundary, the truncation error becomes
\begin{equation}
	\tau_1=\overline{b_1}-(L\phi^e)_1=\left(-\frac{1}{2}\theta_L^2+\theta_L-\frac{1}{2}\right)\phi_1^{''}-\frac{(1-\theta_L)}{3}\phi_1^{'''}\Delta_x+\mathcal{O}(\Delta_x^2),
\end{equation}
which is evidently degraded to a zeroth-order term. A similar derivation for the right boundary yields
\begin{equation}
	\tau_{N-1}=\left(-\frac{1}{2}\theta_R^2+\theta_R-\frac{1}{2}\right)\phi_{N-1}^{''}+\mathcal{O}(\Delta_x).
\end{equation}

Following the same approach as in the linear case, we arrive at an equivalent system of equations relating the solution error $\xi_i$ and the terms $H_{1/2}$ and $H_{N-1/2}$, with the only difference being the updated expressions due to the quadratic boundary treatment:
\begin{equation}
	H_{1/2}=\left(\frac{2-\theta_L}{\theta_L}\xi_1-\frac{1-\theta_L}{1+\theta_L}\xi_2\right)/\Delta_x,
\end{equation}
and
\begin{equation}
	H_{N-1/2}=\left(\frac{2-\theta_R}{-\theta_R}\xi_{N-1}+\frac{1-\theta_R}{1+\theta_R}\xi_{N-2}\right)/\Delta_x.
\end{equation}
 Solving the set of equations leads to
\begin{equation}
	\begin{split}
		H_{N-1/2} &= \Delta_x \left\{ \frac{\theta_L(1+\theta_L)}{2}\tau_1 + \left[(N+\theta_L-2) + \frac{\theta_R(1-\theta_R)}{2}\right]\tau_{N-1} \right. \\
		&\quad \left. - \sum_{k=2}^{N-2}(1-\theta_L-k)\tau_k \right\}/\left(N+\theta_L+\theta_R-2\right),
	\end{split}
\end{equation}
and consequently,
\begin{equation}
	\xi_i=\Delta_x^2\left[(i+\theta_L-1)\frac{H_{N-1/2}}{\Delta_x}-\frac{1}{2}\theta_L(1+\theta_L)\tau_1+\sum_{k=2}^{N-1}(1-\theta_L-k)\tau_k-\sum_{k=i+1}^{N-1}(i-k)\tau_k\right].
\end{equation}

Notably, $H_{N-1/2}$ is now $O(\Delta_x)$—in contrast to $O(\Delta_x^2)$ in absence of RHS inaccuracies\cite{jomaa2005embedded}—due to the degradation of $\tau_1$ and $\tau_{N-1}$. As a result, $\xi_i$ is uniformly second-order in $\Delta_x$, including the boundary nodes $\xi_1$ and $\xi_{N-1}$, which were third-order when the RHS was accurate.

To isolate the left boundary contribution, we set $\theta_R=1$ and retain only the $\tau_1$ term in the above expression, leading to
\begin{equation}
	\xi_i^L=\frac{1}{2}\left(\frac{i}{N}-1\right)\theta_L(1+\theta_L)\tau_1\Delta_x^2\simeq\frac{1}{2}\left(\frac{i}{N}-1\right)\theta_L(1+\theta_L)\left(-\frac{1}{2}\theta_L^2+\theta_L-\frac{1}{2}\right)\phi_1^{''}\Delta_x^2,
\end{equation}
which varies linearly with $i$ and decreases from the left to the right boundary. Further analysis indicates that the left boundary error is maximized when $\theta_L\approx0.3904$. Similarly, for the right boundary contribution we obtain
\begin{equation}
	\xi_i^R=-\frac{i}{2N}\theta_R(1+\theta_R)\tau_{N-1}\Delta_x^2\simeq-\frac{i}{2N}\theta_R(1+\theta_R)\left(-\frac{1}{2}\theta_R^2+\theta_R-\frac{1}{2}\right)\phi_{N-1}^{''}\Delta_x^2,
\end{equation} 
which is also linear in $i$ and reaches its maximum at $\theta_R\approx0.3904$. Again, we observe that the boundary-induced error components $\xi_i^{L}$ and $\xi_i^R$ have signs opposite to those of the corresponding boundary truncation errors $\tau_1$ and $\tau_{N-1}$. This intriguing behavior also emerges in the 2D cases, as confirmed by later numerical experiments.

From these expressions, we conclude that under the influence of inaccurate RHS values, the quadratic scheme exhibits boundary error contributions of order $O(\Delta_x^2)$---qualitatively resembling the behavior of the linear scheme. Meanwhile, the contributions from interior nodes remain $O(\Delta_x^2)$ and can still be approximated by the integral formula~\ref{eqn:inner_contrib}. Thus, the simple decomposition $\xi_i \approx \xi_i^L + \xi_i^R + \xi_i^{in}
$ now holds for the quadratic scheme affected by RHS inaccuracies. In contrast, when the RHS is accurate, the solution error is primarily dominated by the interior contributions\cite{jomaa2005embedded}.

\subsubsection{Overview of 1-D error analysis}
\begin{figure}
	\centering
	\includegraphics[width=0.7\linewidth]{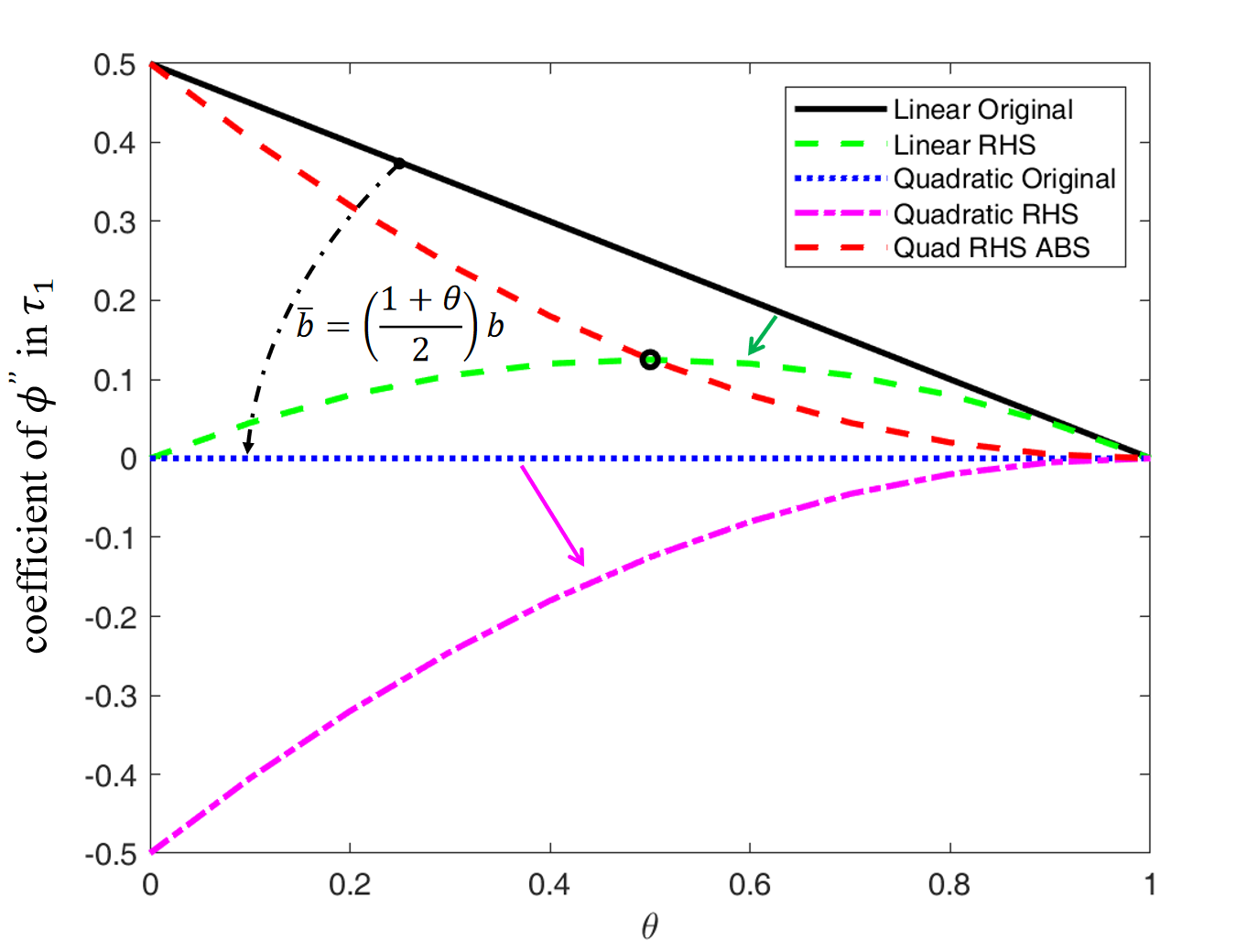}
	\caption{comparison of the leading error terms in $\tau_1$ among four scenarios.}
	\label{fig:tau_compare}
\end{figure}

To summarize the 1-D error analysis, we begin by comparing the leading-order truncation error at the near-boundary node, specifically the coefficient of $\phi_1^{''}$ in $\tau_1$, across four representative scenarios: the linear scheme with and without RHS inaccuracies, and the corresponding quadratic cases. The comparison is illustrated in Figure~\ref{fig:tau_compare}.

For $\theta$ ranging from 0 to 1, the presence of inaccurate RHS values reduces the leading truncation error in $\tau_1$ for the linear scheme, while it increases the error for the quadratic scheme. The plot also indicates a special case in which RHS values are modified as $\overline{b_i} = \left[(1+\theta)/2\right]b_i$ in the linear scheme. This modification effectively transforms the linear scheme to the quadratic accuracy, which can be verified via equations~\ref{eqn:standard}, \ref{eqn:linear}, and \ref{eqn:quad}. Notably, the coefficient for the RHS-affected quadratic case has the opposite sign compared to the linear cases; hence, for comparison purposes, its absolute values are also plotted. The figure shows that for $\theta \in (0, 0.5)$, the RHS-affected linear scheme exhibits a smaller leading error term, whereas for $\theta \in (0.5, 1)$, the RHS-affected quadratic scheme performs better.

\begin{figure}
	\centering
	\includegraphics[width=0.7\linewidth]{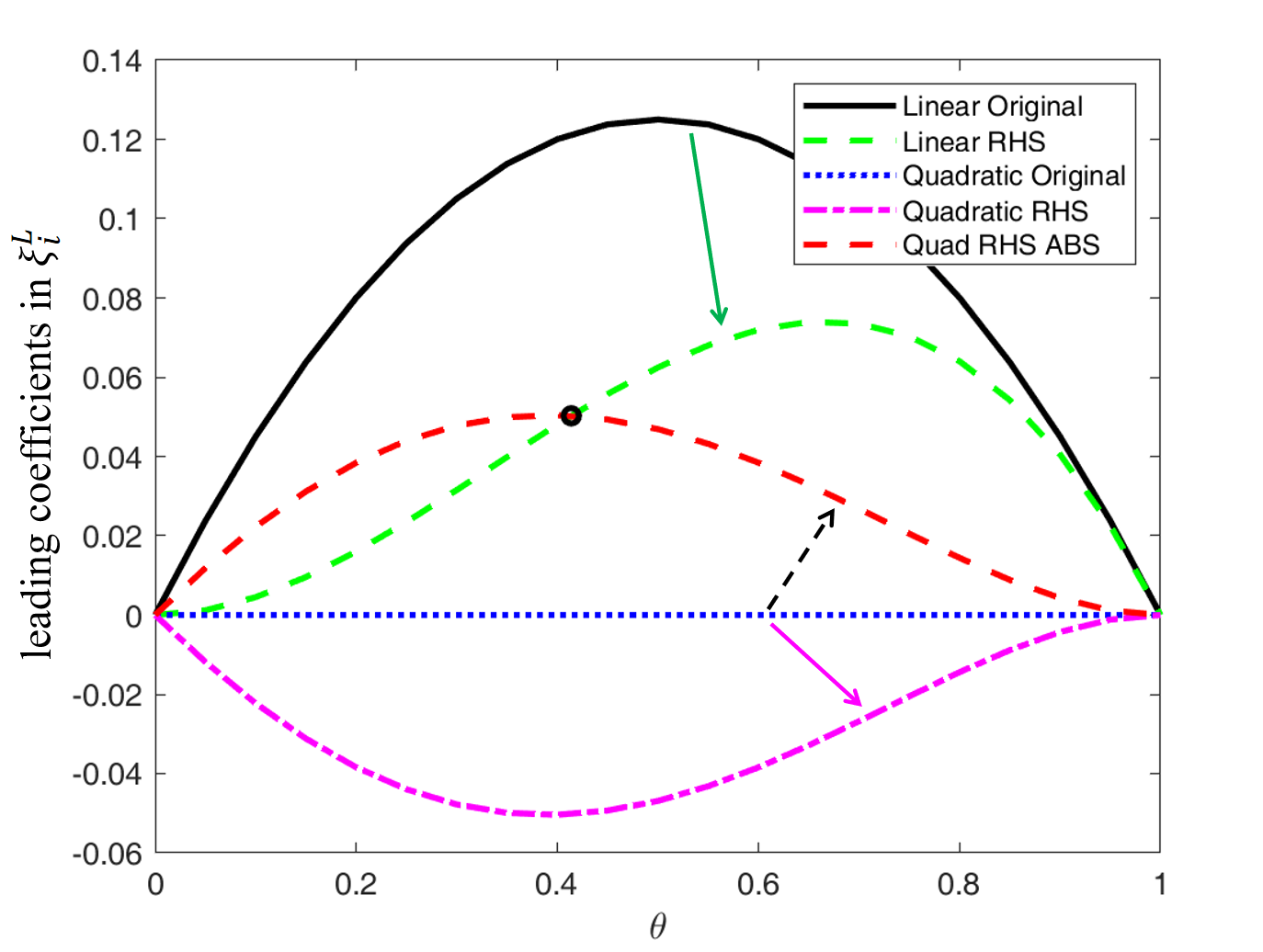}
	\caption{comparison of the leading error terms in $\xi_i^L$ among the four scenarios.}
	\label{fig:xi_compare}
\end{figure}

As established in earlier sections, the magnitudes of $\tau_1$ and $\tau_{N-1}$, in conjunction with $\theta_L$ and $\theta_R$, directly influence the numerical errors $\xi_i$ across the domain through the boundary error components $\xi_i^L$ and $\xi_i^R$. Since the interior node contributions remain the same across all schemes, we focus our comparison on the leading-order boundary error terms---specifically, the coefficients of $(i/N - 1)\phi_1^{''}\Delta_x^2$ in $\xi_i^L$. The results are presented in Figure~\ref{fig:xi_compare}. Note that for the linear scheme unaffected by the RHS inaccuracies, we adopt the $\xi_i^L$ formula from \cite{jomaa2005embedded}; for the unaffected quadratic scheme, since the boundary error is third-order (compared to second-order for the remaining cases), it is approximated as zero.

It's evident that for the linear scheme, boundary errors are reduced across the full range $\theta \in (0, 1)$ when RHS inaccuracies are present. In contrast, the quadratic scheme experiences increased boundary errors under the same effect. Moreover, the boundary error terms for the RHS-affected quadratic scheme and the linear scheme have opposite signs. By comparing the absolute magnitudes of these terms, it is observed that the RHS-affected linear scheme yields smaller boundary errors for $\theta \in (0, 0.414)$, while the RHS-affected quadratic scheme performs better for $\theta \in (0.414, 1)$. Interestingly, this critical $\theta$ value is slightly different from that identified in the $\tau_1$ comparison.

In conclusion, the 1-D error analysis clearly demonstrates that both the choice of boundary treatment scheme and the RHS inaccuracies significantly influence the numerical solution accuracy for the Poisson equation. These effects manifest primarily by modifying the magnitude and order of the boundary truncation errors $\tau_1$ and $\tau_{N-1}$, which in turn shape the boundary error components $\xi_i^L$ and $\xi_i^R$. Specifically, RHS inaccuracies tend to improve the linear scheme by reducing the leading truncation errors at boundaries, while they degrade the accuracy of the quadratic scheme by introducing zeroth-order error components to the truncation errors at boundaries.

\subsection{2-D error analysis}
\begin{figure}
	\centering
	\includegraphics[width=0.6\linewidth]{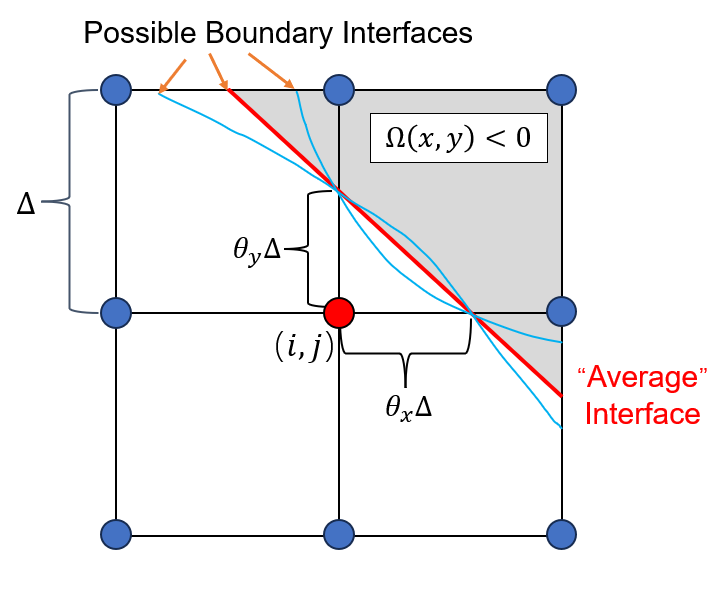}
	\caption{Illustration of a typical near-boundary node in 2D problems that is imperfectly positioned along both axial directions.}
	\label{fig:2D_node}
\end{figure}

For 2D (or 3D) problems, explicit numerical error expressions are difficult to obtain. Instead, we focus on analyzing how inaccurate RHS values affect the truncation errors at near-boundary nodes. Consider a typical near-boundary node $(i,j)$ illustrated in figure~\ref{fig:2D_node}, where $\theta_x$ and $\theta_y$ represent the relative distances to the boundary along the $x$- and $y$-directions, respectively. For a given pair $(\theta_x,\theta_y)$, there exist infinitely many possible interface shapes. To enable tractable analysis, we assume that made that when the mesh spacing $\Delta$ is sufficiently small, the ``average" interface can be approximated as a straight line uniquely determined by $\theta_x$ and $\theta_y$. 

In the absence of RHS inaccuracies, the truncation errors at node $(i,j)$ is defined as
\begin{equation}
	\tau_{i,j}=b_{i,j}-(L\phi^e)_{i,j},
\end{equation}
where the Laplacian operator $L$ acts on both $x$ and $y$ directions. If $(i,j)$ is an interior node (i.e., $\theta_x=\theta_y=1$), we can expand $\phi(x,y)$ around the spatial coordinate $(x_i,y_j)$ to obtain 
\begin{equation}
	\tau_{i,j}=-\frac{1}{12}\left[\frac{\partial^4\phi}{\partial x^4}\bigg|_{(x_i, y_j)}+\frac{\partial^4\phi}{\partial y^4}\bigg|_{(x_i, y_j)}\right]\Delta^2+\mathcal{O}(\Delta^4),
\end{equation}
which is second-order in $\Delta$. However, when node $(i,j)$ is a near-boundary node, the chosen boundary treatment significantly alters the truncation error. Under the linear scheme, the truncation error becomes
  \begin{equation}
  	\label{eqn:tau_linear}
  	\tau_{linear}=\frac{1-\theta_x}{2}\frac{\partial^2\phi}{\partial x^2}\bigg|_{(x_i, y_j)}+\frac{1-\theta_y}{2}\frac{\partial^2\phi}{\partial y^2}\bigg|_{(x_i, y_j)}+\mathcal{O}(\Delta),
  \end{equation}
  which is zeroth-order with respect to $\Delta$. If adopting the quadratic scheme, we have
    \begin{equation}
    	\label{eqn:tau_quad}
  	\tau_{quad}=\left[\frac{1-\theta_x}{3}\frac{\partial^3\phi}{\partial x^3}\bigg|_{(x_i, y_j)}+\frac{1-\theta_y}{3}\frac{\partial^3\phi}{\partial y^3}\bigg|_{(x_i, y_j)}\right]\Delta+\mathcal{O}(\Delta^2),
  \end{equation}
  which is first-order.
  
Now, consider the impact of inaccurate RHS values under the assumption of locally uniform RHS distribution. For the linear scheme, the truncation error becomes:
  \begin{eqnarray}
  	\label{eqn:tau_linear_RHS}
  	\tau_{linear}^{RHS}&=&\overline{b_{i,j}}-(L\phi^e)_{i,j} = b_{i,j}-(L\phi^e)_{i,j}+\delta_{i,j} b_{i,j}-b_{i,j} \nonumber\\
  	&=& \left(\frac{1-\theta_x}{2}+\delta_{i,j}-1\right)\frac{\partial^2\phi}{\partial x^2}\bigg|_{(x_i, y_j)}+\left(\frac{1-\theta_y}{2}+\delta_{i,j}-1\right)\frac{\partial^2\phi}{\partial y^2}\bigg|_{(x_i, y_j)},
  \end{eqnarray}
 where only leading-order terms are retained. Similarly, for the quadratic scheme we have
    \begin{eqnarray}
    	\label{eqn:tau_quad_RHS}
  	\tau_{quad}^{RHS}&=&\overline{b_{i,j}}-(L\phi^e)_{i,j} \nonumber\\&=& \left(\delta_{i,j}-1\right)\left[\frac{\partial^2\phi}{\partial x^2}\bigg|_{(x_i, y_j)}+\frac{\partial^2\phi}{\partial y^2}\bigg|_{(x_i, y_j)}\right],
  \end{eqnarray}
which now reduces to zeroth-order. From these expressions, it is evident that the quadratic scheme experiences a loss of accuracy due to degeneration of its boundary truncation error when inaccurate RHS values are present. In contrast, the effect on the linear scheme is more nuanced and depends on the value of $\delta_{i,j}$, which warrants further numerical evaluation.
  
\begin{figure}[htbp]
  	\centering
  	\begin{subfigure}[b]{0.32\textwidth}
  		\includegraphics[width=\linewidth]{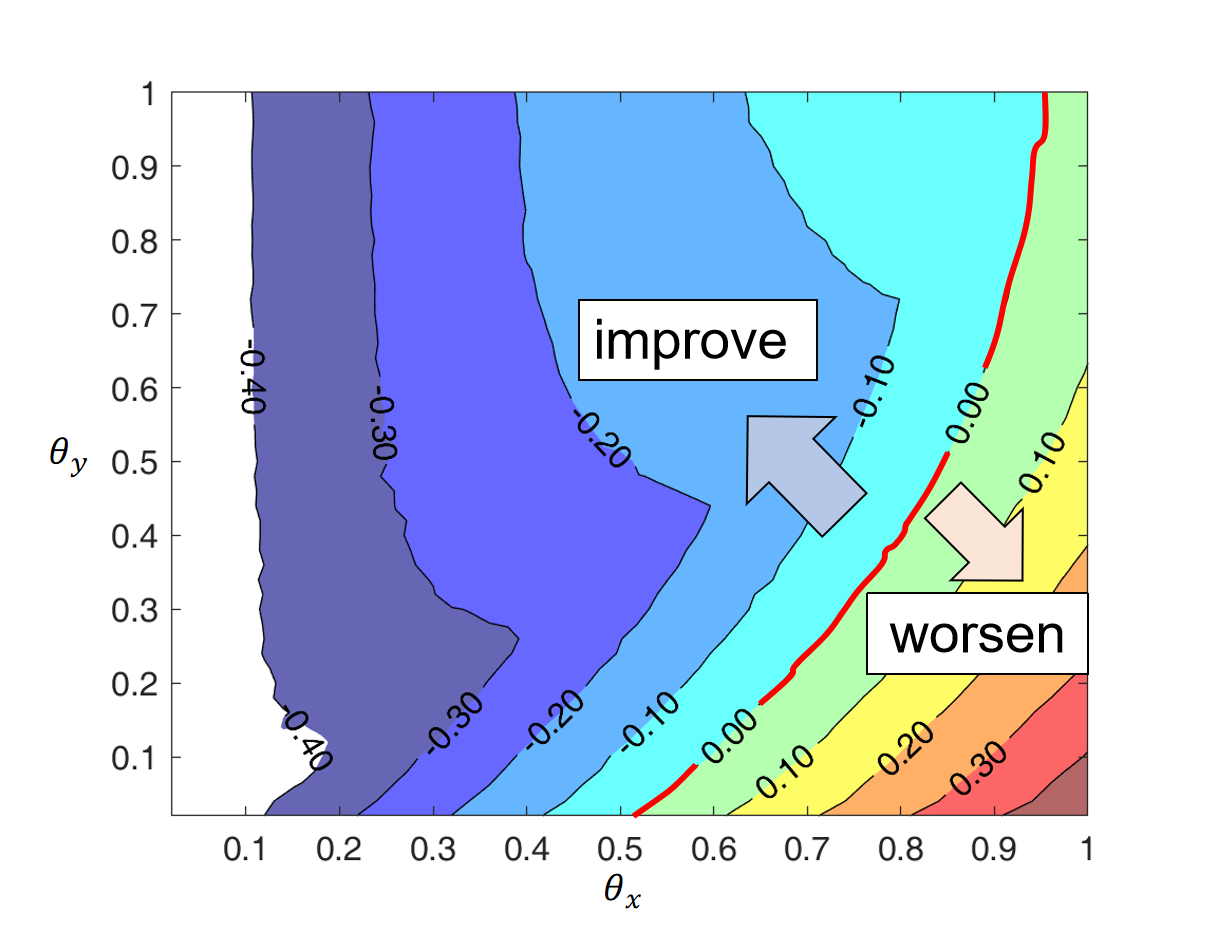}
  		\caption{}
  		\label{fig:linear_2D_tau_a}
  	\end{subfigure}
  	\hfill 
  	\begin{subfigure}[b]{0.32\textwidth}
  		\includegraphics[width=\linewidth]{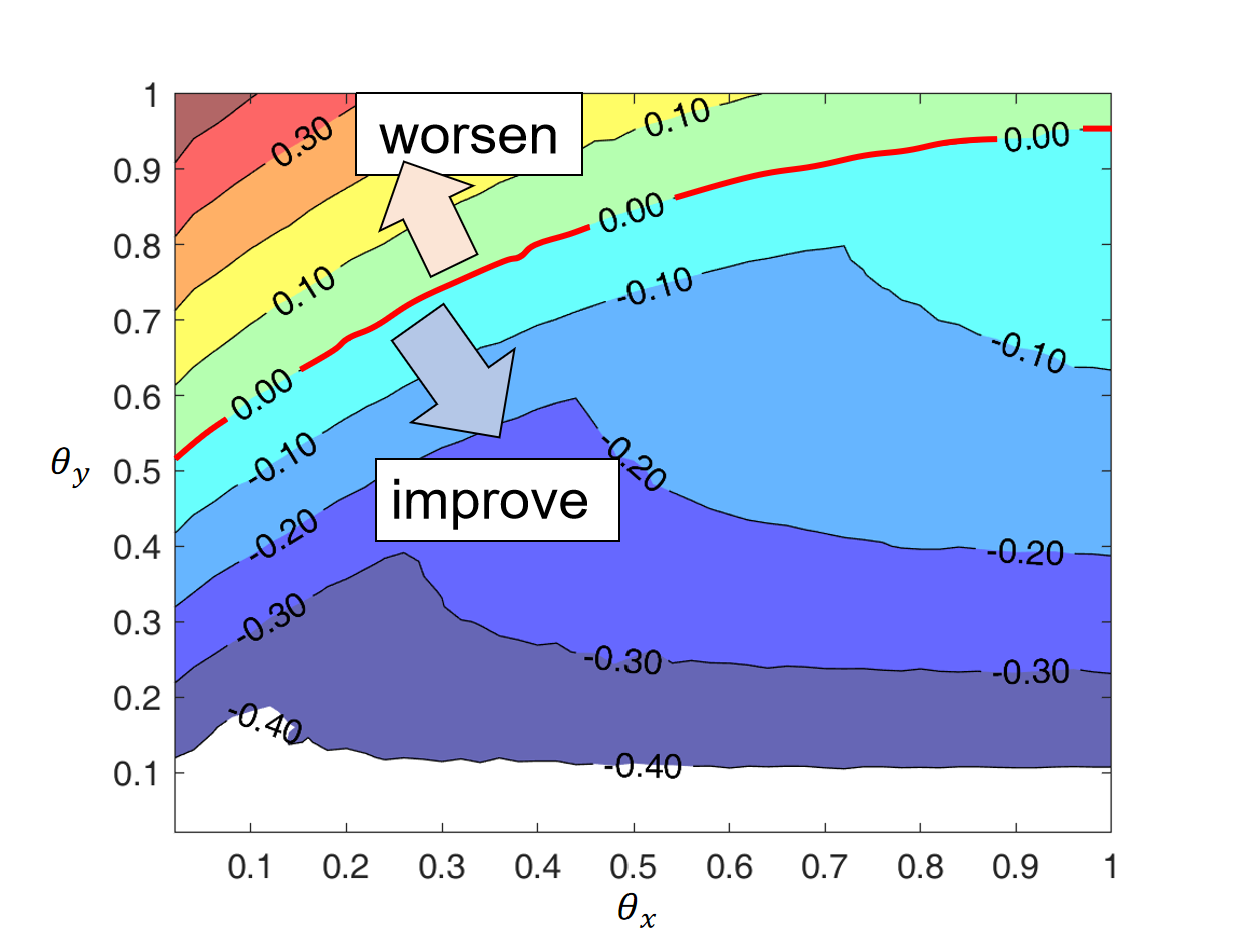}
  		\caption{}
  		\label{fig:linear_2D_tau_b}
  	\end{subfigure}
  	\hfill
  	\begin{subfigure}[b]{0.32\textwidth}
  		\includegraphics[width=\linewidth]{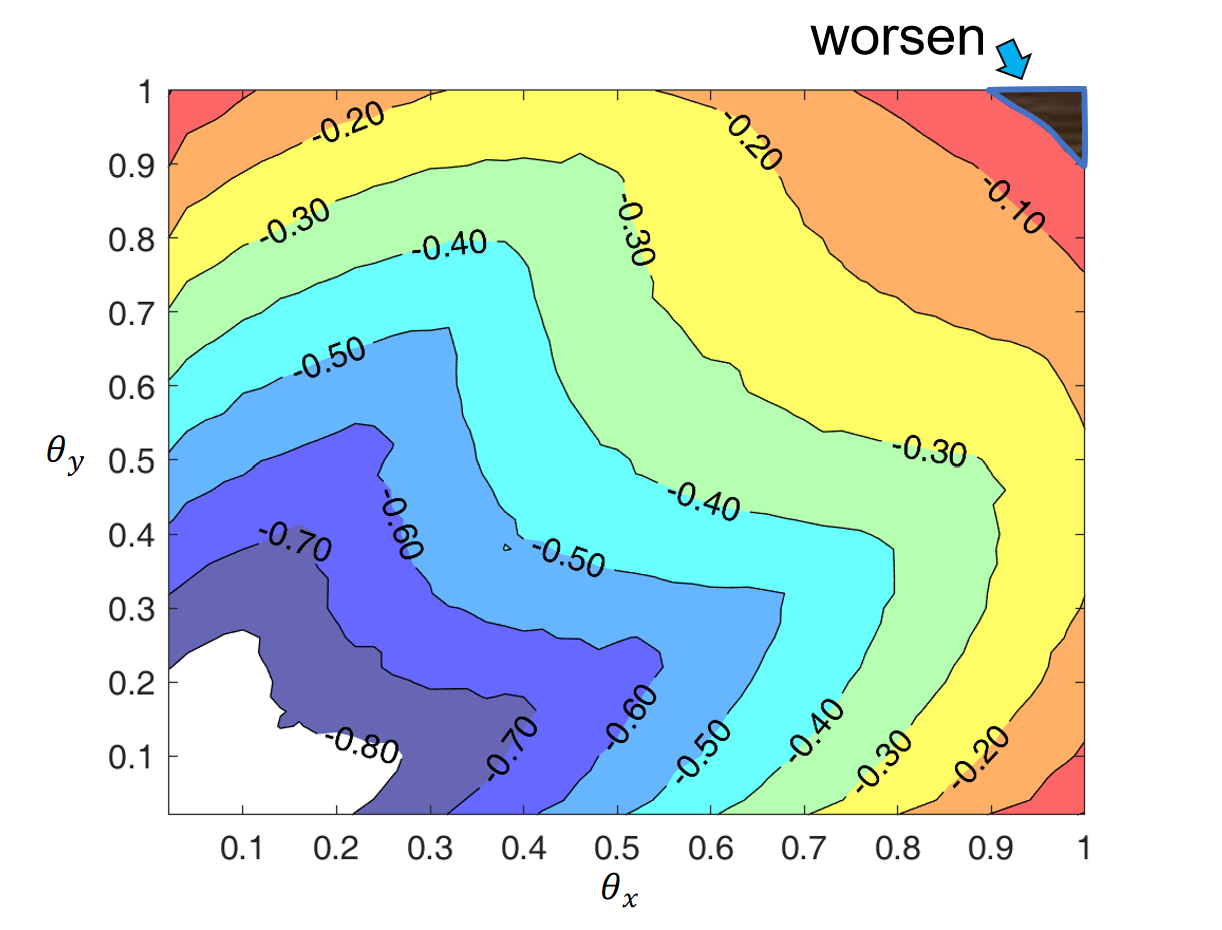}
  		\caption{}
  		\label{fig:linear_2D_tau_c}
  	\end{subfigure}
  	\caption{The contour of the linear scheme's leading error term coefficients comparison with/without inaccurate RHS effect across all range of $\theta_x$ and $\theta_y$. $(a)$ compares the coefficients of the  $\partial^2\phi/\partial x^2$ terms by subtraction, where negative values means RHS-affected linear scheme has lower truncation errors; $(b)$ compares the $\partial^2\phi/\partial y^2$ terms; $(c)$ combines the results of $(a)$ and $(b)$, assuming $\partial^2\phi/\partial x^2\simeq  \partial^2\phi/\partial y^2$.}
  	\label{fig:linear_2D_tau_comparison}
  \end{figure}

To this end, we numerically compute the approximate $\overline{\delta_{i,j}}$ across the full range of $\theta_x$ and $\theta_y$ using the method described in equation\eqref{eqn:delta}. We then examine how the coefficients of the $\partial^2 \phi / \partial x^2$ and $\partial^2 \phi / \partial y^2$ terms in the truncation error are affected in the linear scheme. The results are shown in Figure~\ref{fig:linear_2D_tau_comparison}. It can be observed from both Figures~\ref{fig:linear_2D_tau_a} and \ref{fig:linear_2D_tau_b} that, for a given $(\theta_x, \theta_y)$ pair, RHS inaccuracies tend to improve the linear scheme by reducing the coefficients of the leading-order error terms. Further assuming $\partial^2 \phi / \partial x^2 \simeq \partial^2 \phi / \partial y^2$, Figure~\ref{fig:linear_2D_tau_c} shows that the RHS-affected linear scheme consistently reduces the boundary truncation error, except for a narrow band where both $\theta_x$ and $\theta_y$ are close to 1.
  
\begin{figure}[htbp]
  	\centering
  	\begin{subfigure}[b]{0.32\textwidth}
  		\includegraphics[width=\linewidth]{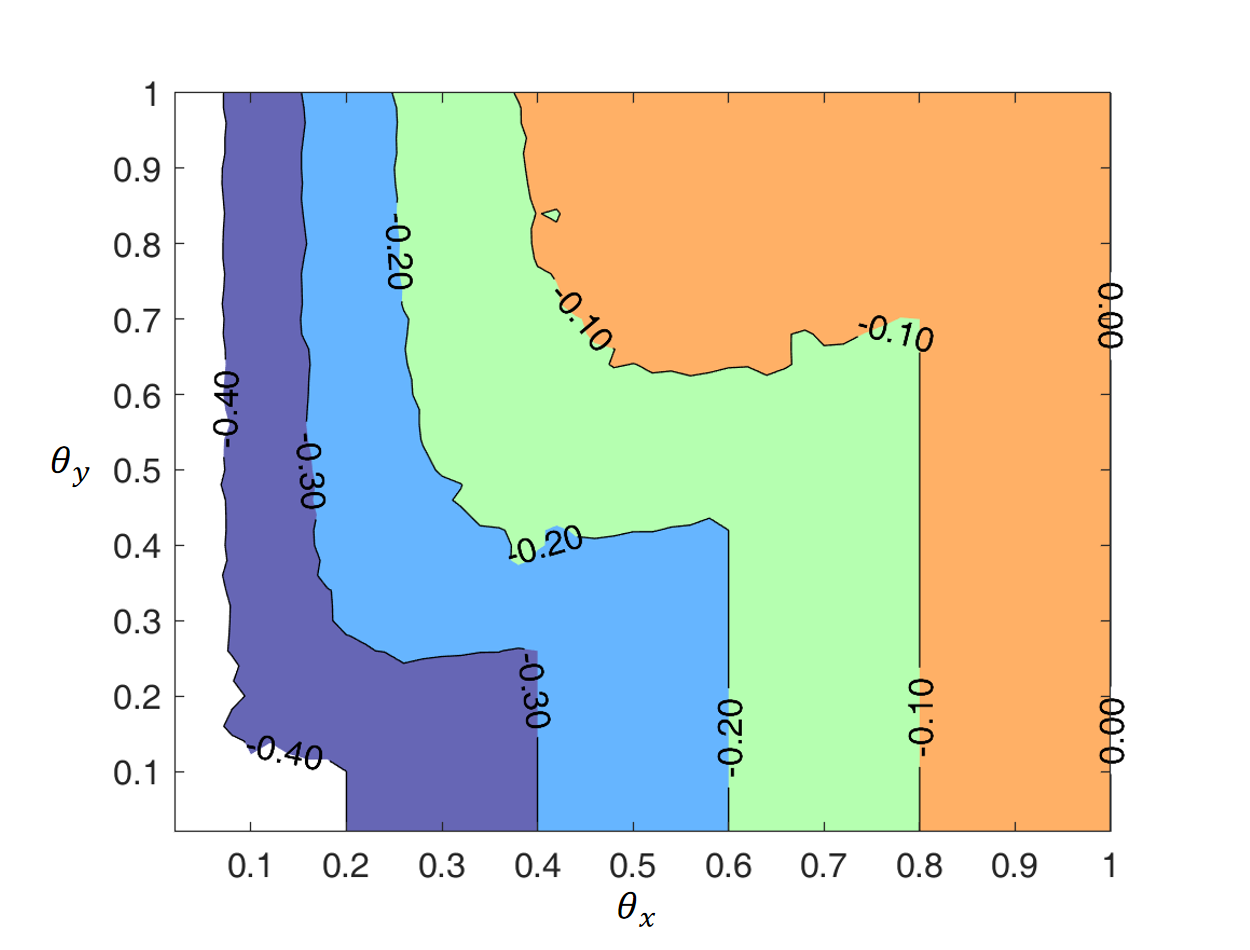}
  		\caption{}
  		\label{fig:LvsQ_2D_tau_a}
  	\end{subfigure}
  	\hfill 
  	\begin{subfigure}[b]{0.32\textwidth}
  		\includegraphics[width=\linewidth]{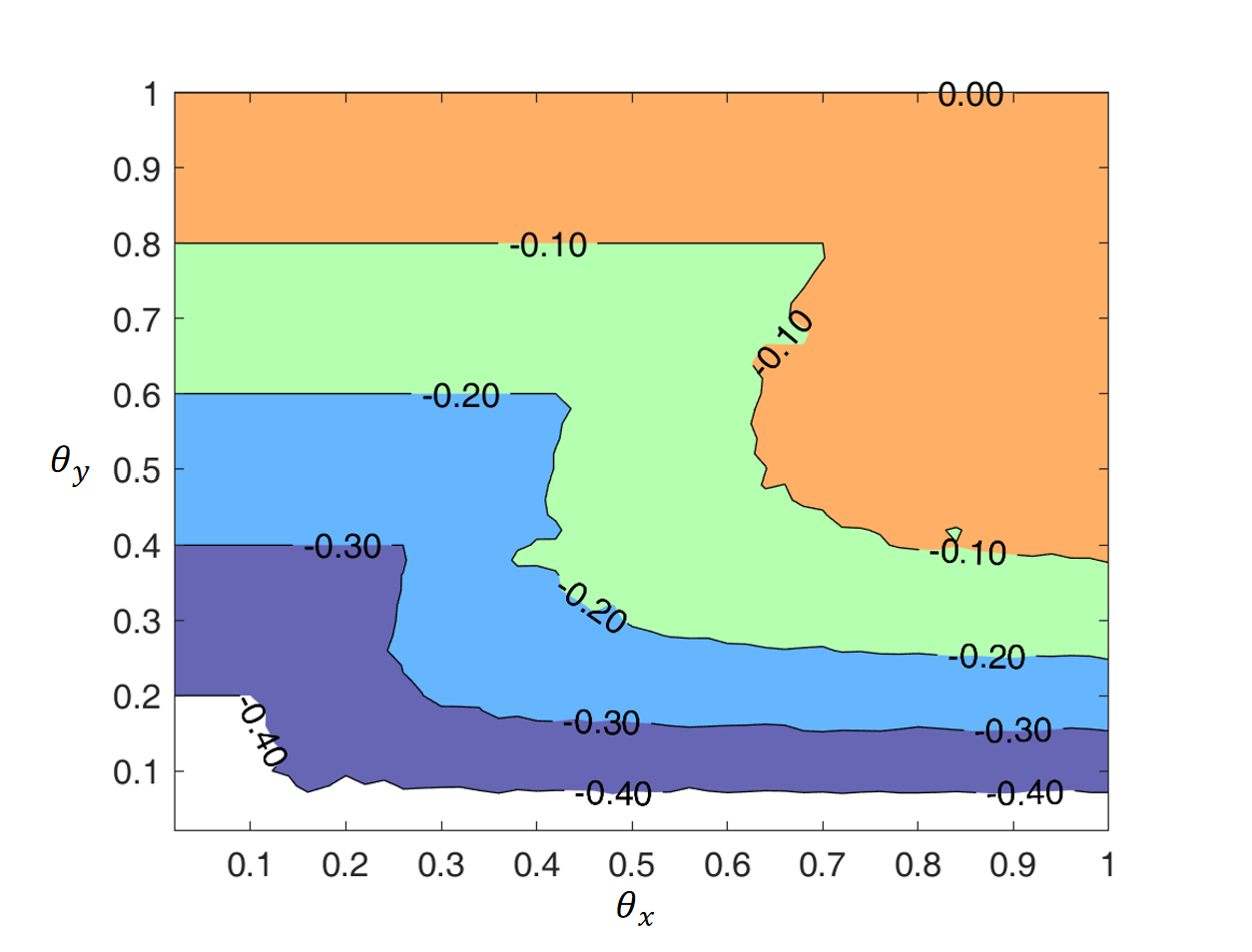}
  		\caption{}
  		\label{fig:LvsQ_2D_tau_b}
  	\end{subfigure}
  	\hfill
  	\begin{subfigure}[b]{0.32\textwidth}
  		\includegraphics[width=\linewidth]{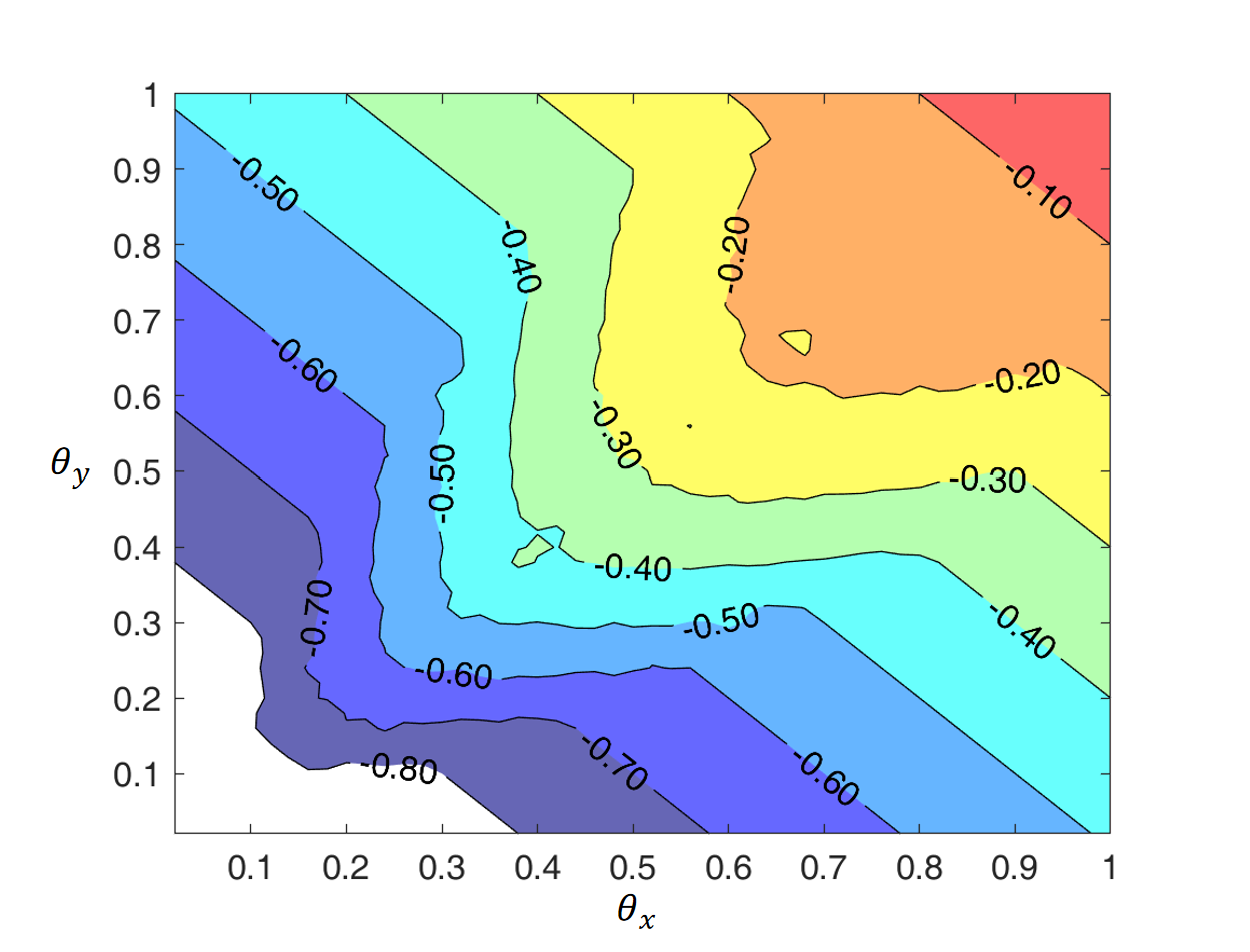}
  		\caption{}
  		\label{fig:LvsQ_2D_tau_c}
  	\end{subfigure}
  	\caption{The contour of leading error term coefficients comparison between the linear and the quadratic scheme under inaccurate RHS effect across all range of $\theta_x$ and $\theta_y$. $(a)$ compares the coefficients of the  $\partial^2\phi/\partial x^2$ terms by subtraction, where negative values means RHS-affected linear scheme has lower truncation errors; $(b)$ compares the $\partial^2\phi/\partial y^2$ terms; $(c)$ combines the results of $(a)$ and $(b)$, assuming $\partial^2\phi/\partial x^2\simeq  \partial^2\phi/\partial y^2$.}
  	\label{fig:LvsQ_2D_tau_comparison}
  \end{figure}

We further compare the boundary truncation errors between the linear and quadratic schemes, both under the influence of inaccurate RHS values. As shown in Figure~\ref{fig:LvsQ_2D_tau_comparison}, the RHS-affected linear scheme exhibits uniformly smaller truncation errors than the quadratic scheme across the entire $(\theta_x, \theta_y)$ domain. However, further numerical experiments are necessary to determine whether this advantage in boundary truncation errors translates to overall improvements in the solution accuracy of the Poisson equation, particularly in scenarios where RHS values are not uniformly distributed.

\section{Test cases and results}
\label{sec:test case}
Numerical tests were conducted on a series of 1D, 2D, and 3D examples with predefined analytical solutions $\phi(\vec{x})$. For each test case, the computational domain was embedded within a uniformly spaced Cartesian mesh, and both the linear and quadratic boundary treatment schemes were employed to discretize the Poisson equation, resulting in a linear system of the form $\mathbf{A}\vec{\phi} = \vec{b}$. The system was solved using a direct solver with row pivoting to reduce rounding errors, rather than iterative solvers. After obtaining the numerical solution, the error vector was computed as $\vec{\xi} = \vec{\phi} - \vec{\phi}^e$.

While various boundary geometries and analytical forms of $\phi(\vec{x})$ were also tested, the results consistently supported the same conclusions. Therefore, for brevity and clarity, we present only three representative examples in this section. The test cases are described as follows.

\begin{enumerate}
    \item 1D example. Poisson equation is to be solved in a 1D domain $[-0.5,0.5]$ with the exact solution given by $\phi(x)=4x^2\sin{2\pi x}$, which determines the RHS values by $b(x)=d^2\phi/dx^2=8(1-2\pi^2x^2)\sin{2\pi x}+32\pi x\cos{2\pi x}$. Dirichlet boundaries are defined such that $|x|\geq0.313$ denotes the outside region.
    \item 2D example. 2D Poisson equation is defined within a square area $[-1,1]\times[-1,1]$ with an exact solution $\phi(x,y)=\left[ (x+2)^2+(y-2)^2\right]^{-1}$, which has the corresponding RHS values of $b(x,y)=\nabla^2\phi=4/(8+4x+x^2-4y+y^2)^2$. A closed Dirichlet boundary interface is defined by a pair of parameterized equations $x(t)$ and $y(t)$, where $x(t)=0.02\sqrt{5}+[0.5+0.2\sin{5t}]\cos{t}$, and $y(t)=0.02\sqrt{5}+[0.5+0.2\sin{5t}]\sin{t}$. This Dirichlet interface is shaped like a starfish, whose inside contains all the inner grid nodes. In order to decide whether a grid node $(i,j)$ is inside the boundary and to calculate the $\theta_x$ and $\theta_y$ if it's a near-boundary node, it is recommended that the calculations are done in a shifted polar coordinate system originated at $(0.02\sqrt{5},0.02\sqrt{5})$ , then any point that is within the radius of $(0.5+0.2\sin{5t})$ is inside the boundary.
    \item 3D example. The 3D Poisson equation is solved within a cubic domain $[0,1]^3$, where a spherical boundary decides the inner nodes by $\Omega<0$, where $\Omega(x,y,z)=\sqrt{(x-0.5)^2+(y-0.5)^2+(z-0.5)^2}-0.3$. The exact solution is given by $\phi(x,y,z)=\exp{(-x^2-y^2-z^2)}$, from which the RHS values are computed as $b(x,y,z)=\exp{(-x^2-y^2-z^2)}(4x^2+4y^2+4z^2-6)$.
\end{enumerate}

\subsection{1-D numerical test results}
We numerically solved the 1-D Poisson example under four scenarios: using accurate RHS values $b_i=b(x_i)$ and inaccurate RHS values $\overline{b_i}=\overline{\delta_i}b_i$ (with $\overline{\delta_i}$ computed from equation\eqref{eqn:1D_delta}), each under the linear and quadratic treatments. For each case, we also computed the contributions to the numerical error $\xi_i$ from three components: left boundary ($\xi_i^L$), right boundary ($\xi_i^R$), and inner nodes ($\xi_i^{\text{in}}$). Notably, the inner node contribution $\xi_i^{\text{in}}$ remains unchanged across all cases. These decomposed components were used to validate the approximate relationship $\xi_i \approx \xi_i^L + \xi_i^R + \xi_i^{\text{in}}$, and were plotted together for comparison.

\begin{figure}
	\centering
	\includegraphics[width=0.7\linewidth]{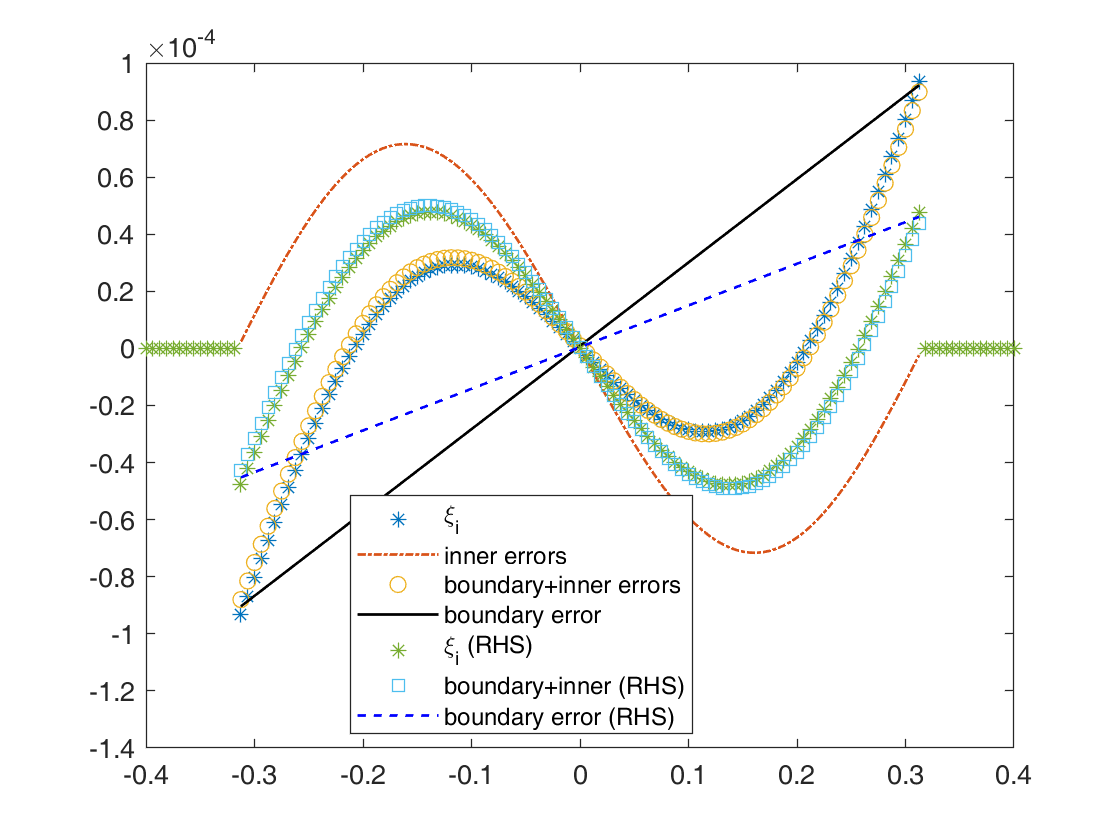}
	\caption{The numerical errors from solving the 1-D Poisson equation with and without RHS effect by the linear scheme. A total number of 161 nodes were set evenly within the domain $[-0.5,0.5]$, the Dirichlet boundaries were set at $\pm0.3156$ to give $\theta_L,\theta_R$=0.5.}
	\label{fig:linear_vs_linearRHS_1D}
\end{figure}

Figure~\ref{fig:linear_vs_linearRHS_1D} compares the error distributions obtained from the linear scheme with and without RHS inaccuracies. It can be seen that the numerical error at each node is well approximated by the sum of boundary and inner  contributions. By introducing RHS effect, the boundary contributions (represented by two linear lines) were reduced at all nodes, but more evident at the two boundary ends. This reduction directly translates into a lower numerical error at the boundary nodes. In general, the reduced boundary contributions lead to an overall decrease in $\xi_i$ at each node. However, it is worth noting that the signs of $\xi_i^L$ and $\xi_i^R$ can be opposite to that of $\xi_i^{\text{in}}$, as is the case here, slightly complicating the interpretation of error behavior.

\begin{figure}
	\centering
	\includegraphics[width=0.7\linewidth]{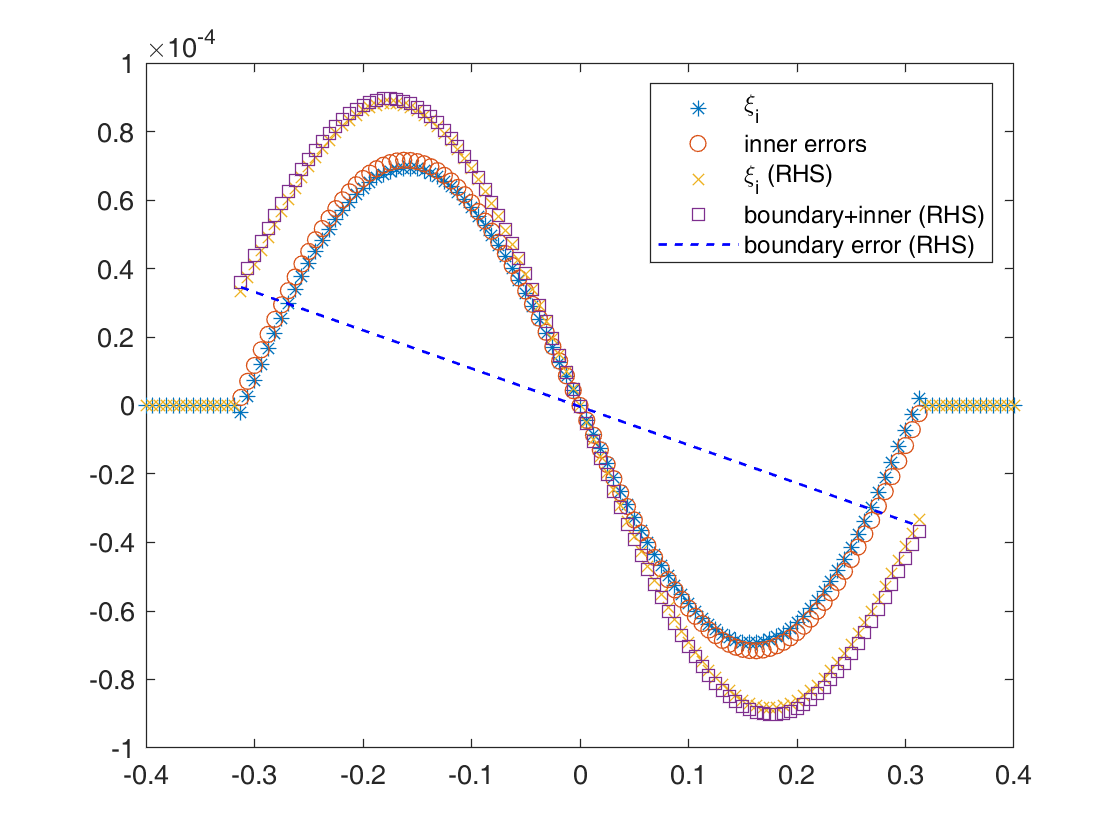}
	\caption{The numerical errors from solving the 1-D Poisson equation with and without RHS effect by the quadratic scheme. A total number of 161 nodes were set evenly within the domain $[-0.5,0.5]$, the Dirichlet boundaries were set at $\pm0.3156$ to give $\theta_L,\theta_R$=0.5, the boundary error contributions only plotted for the RHS case.}
	\label{fig:quad_vs_quadRHS_1D}
	\end{figure}

Figure~\ref{fig:quad_vs_quadRHS_1D} shows a similar comparison for the quadratic scheme. In the RHS-affected case, numerical errors again align with the sum of boundary and inner contributions. Without RHS inaccuracies, the errors are dominated by the inner node contributions alone. The inclusion of inaccurate RHS introduces additional boundary error components, which vary approximately linearly from the left to the right boundary. Interestingly, the signs of the boundary error contributions in this case are opposite to those in the linear scheme, in agreement with the earlier 1-D error analysis. As a result, these boundary terms constructively add to the inner node contributions, amplifying the total numerical error $\xi_i$ throughout the domain.

\begin{figure}
	\centering
	\includegraphics[width=0.7\linewidth]{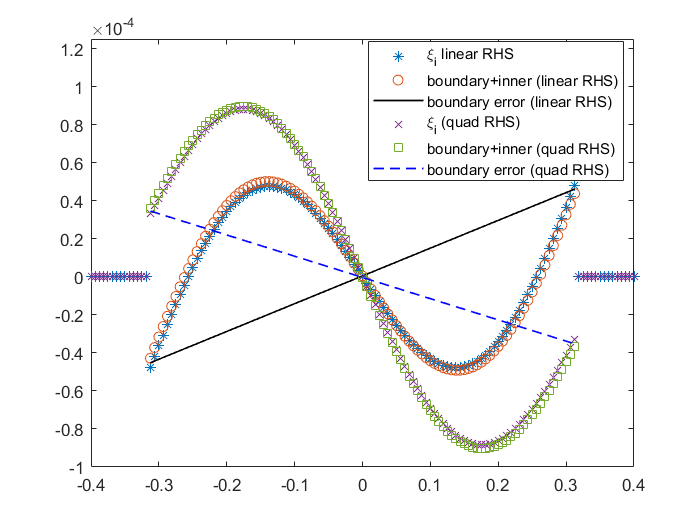}
	\caption{Comparison between the errors from the RHS-affected linear and quadratic scheme solutions. A total number of 161 nodes were set evenly within the domain $[-0.5,0.5]$, the Dirichlet boundaries were set at $\pm0.3156$ to give $\theta_L,\theta_R$=0.5.}
	\label{fig:linearRHS_vs_quadRHS_1D}
\end{figure}

Figure~\ref{fig:linearRHS_vs_quadRHS_1D} provides a direct comparison of the numerical errors obtained using the linear and quadratic schemes under inaccurate RHS values. Both methods exhibit qualitatively similar behavior: numerical errors arise from additive contributions of boundary and interior nodes. The boundary contributions, $\xi_i^L$ and $\xi_i^R$, are most pronounced near the two boundary ends, where the impact of RHS inaccuracies is most significant. Notably, these boundary-induced errors have opposite signs in the linear and quadratic schemes, consistent with the trends observed in the 1-D error analysis. 

\begin{figure}
	\centering
	\includegraphics[width=0.7\linewidth]{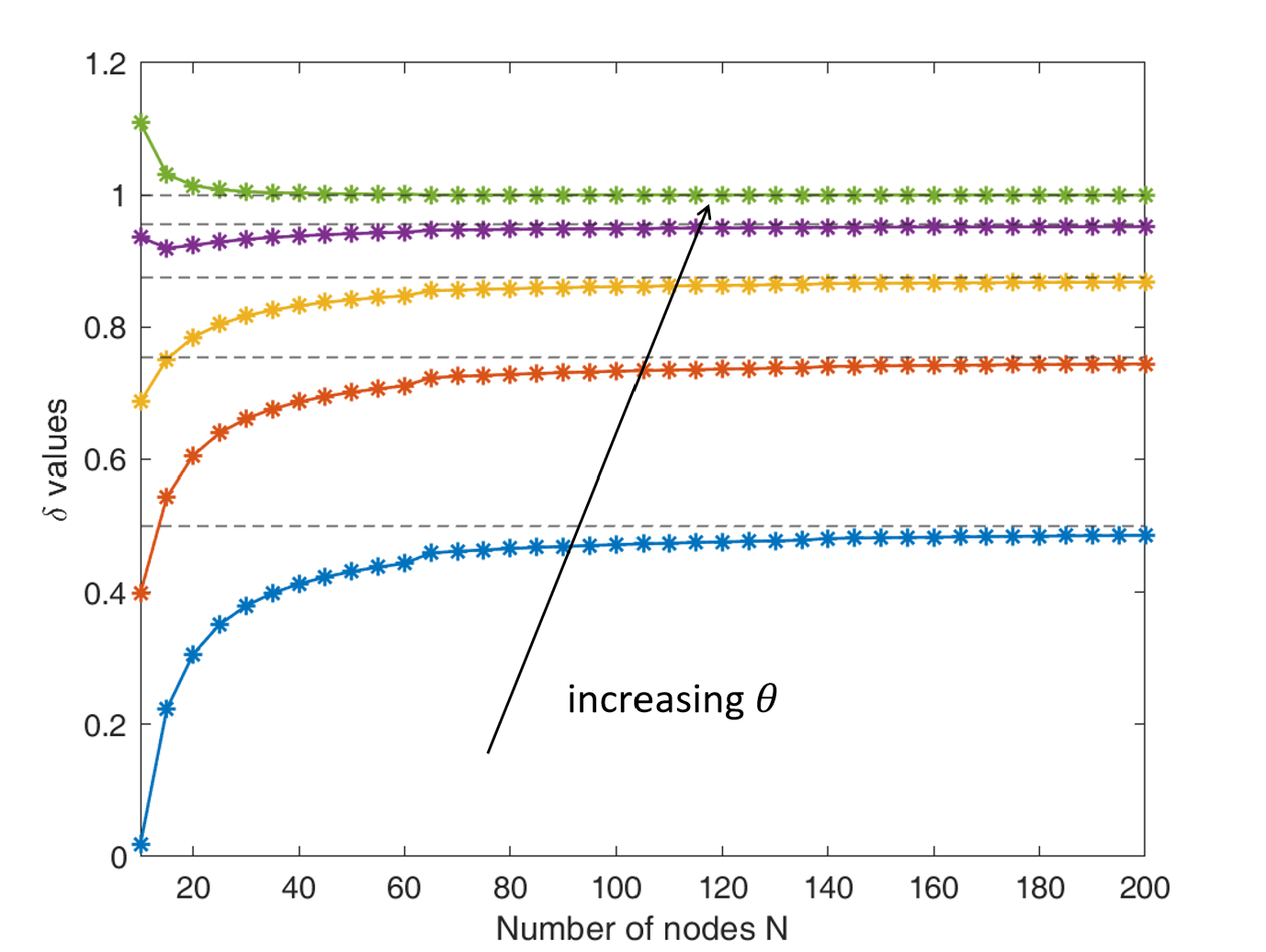}
	\caption{Comparison between $\delta_i$ values calculated by equation\eqref{eqn:1D_delta} and those from more accurate integration equation\eqref{eqn:1D_delta_int}, selected $\theta$ values are 0.0, 0.3, 0.5, 0.7, and 1.0. The horizontal dashed lines represent the accurate $\delta_i$ values calculated by the integral equation, whereas the solid lines represent those calculated by the simpler equation\eqref{eqn:1D_delta}.}
	\label{fig:1D_delta_evolution}
\end{figure}

Overall, the 1-D numerical tests confirmed the findings from the 1-D error analysis in previous sections. We emphasize that for consistency with earlier analysis, the approximated $\overline{\delta_i}$ values used in the tests were computed using equation\eqref{eqn:1D_delta}, rather than the more accurate integral form of equation\eqref{eqn:1D_delta_int}. This simplification has a negligible effect on the results when $\Delta x$ is sufficiently small, as the two expressions yield nearly identical $\delta_i$ values under such conditions. This agreement is illustrated in Figure~\ref{fig:1D_delta_evolution}, which plots both forms of $\delta_i$ for various choices of $\theta$ and $\Delta x$. The close match between the two confirms the validity of using the simpler form in both the numerical experiments and theoretical analysis.

\subsection{2-D/3-D numerical tests and discussions}
For the 2-D and 3-D tests, we followed the same four scenarios as in the 1-D case, involving linear and quadratic schemes with and without inaccurate RHS values. For the RHS-affected cases, the modified RHS values $\overline{b_{i,j}}$ were computed using equation\eqref{eqn:rhs_calc} at all near-boundary nodes. The corresponding $\overline{\delta_{i,j}}$ values were obtained as $\overline{b_{i,j}}/b_{i,j}$, where the accurate RHS values $b_{i,j}$ were evaluated from the prescribed $b(x,y)$ expression. Due to the lack of explicit analytical expressions for boundary and interior error contributions in higher dimensions, we began by analyzing the truncation errors at near-boundary nodes and compared them with the resulting numerical errors, aiming to identify qualitative correlations. For the 2-D case, the full numerical error field was visualized, while for 3-D, representative numerical results were summarized in tabular form. Lastly, for RHS-affected cases, we also tested the effectiveness of using approximated $\overline{\delta_{i,j}}$ values to counteract the impact of RHS inaccuracies.
\subsubsection{Boundary $\tau$ and $\xi$ for the 2-D case}
\begin{figure}
	\centering
	\includegraphics[width=0.7\linewidth]{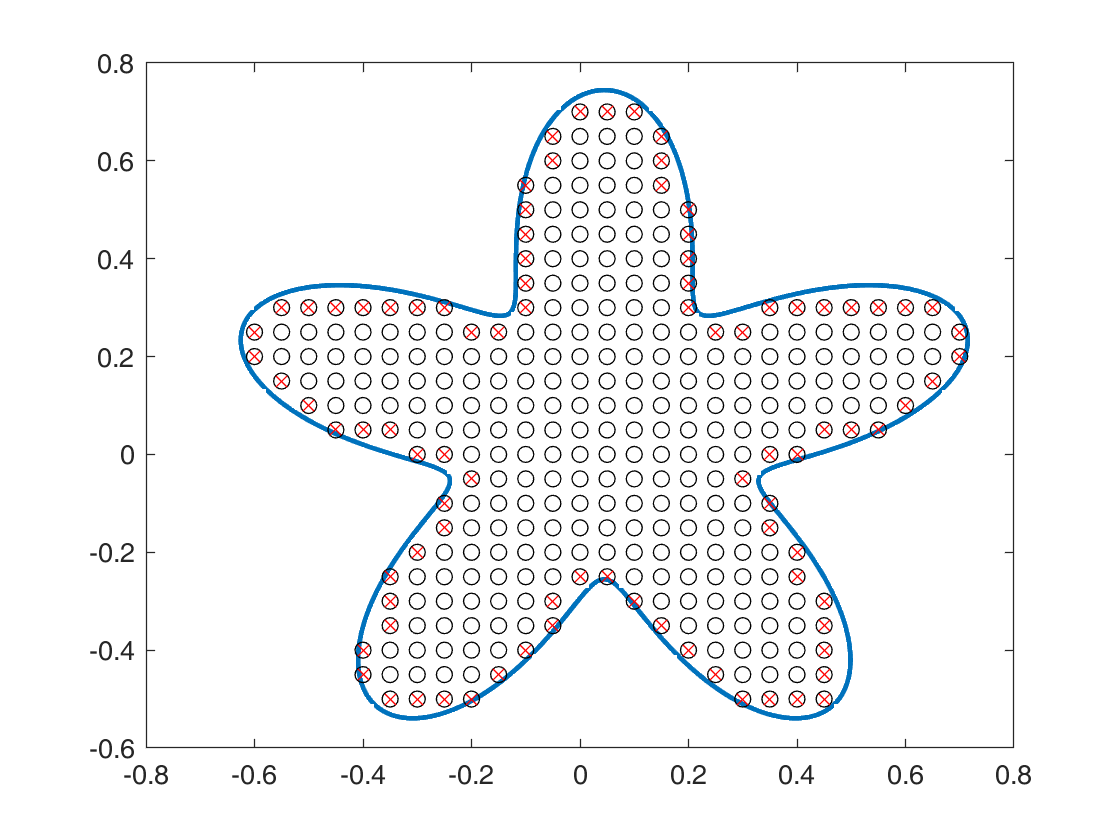}
	\caption{The computational domain for the 2-D case with a total $41\times41$ uniformly spaced nodes, the open circles with a cross indicate the near-boundary nodes.}
	\label{fig:41x41}
\end{figure}

\begin{figure}
	\centering
	\includegraphics[width=0.7\linewidth]{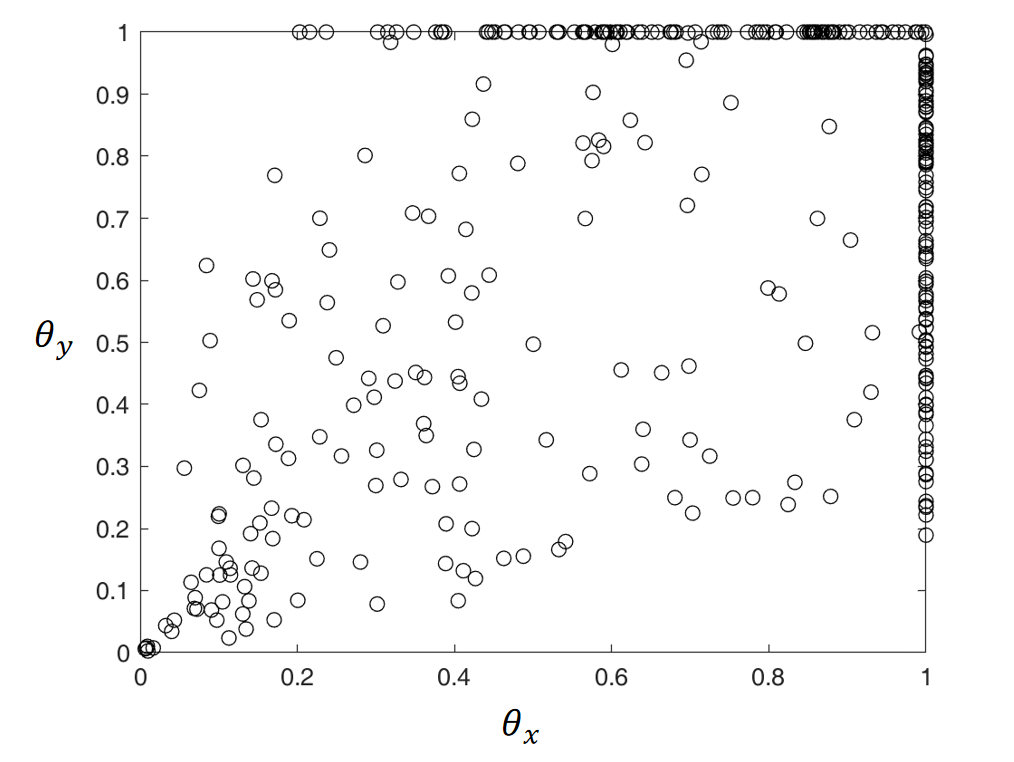}
	\caption{The $\theta_x$ and $\theta_y$ values associated with each near-boundary node. Each open circle corresponds to a particular near-boundary node; a total number of $151\times151$ nodes were used on the computational domain. }
	\label{fig:151x151_theta_map}
\end{figure}

Figure\ref{fig:41x41} illustrates the computational domain for the 2-D case. The solid line denotes the Dirichlet boundary, the open circles indicate the interior nodes, and the circles marked with a cross represent the near-boundary nodes. A total of $41 \times 41$ nodes are used here for demonstration purposes; however, much finer grids were employed in the actual numerical solution of the 2-D Poisson equation. Each near-boundary node is associated with a pair of $\theta_x$ and $\theta_y$ values. Figure~\ref{fig:151x151_theta_map} displays all the $\theta$ pairs when a $151 \times 151$ grid was employed. It is not surprising that the distribution of $\theta$ values is non-uniform, as it strongly depends on the shape of the Dirichlet boundary. Notably, a considerable number of near-boundary nodes have either $\theta_x = 1$ or $\theta_y = 1$.

We first investigated the RHS effect on the linear scheme by examining the truncation errors at the near-boundary nodes. The exact solution $\phi(x,y)$ (and thus its second partial derivatives) enables us to compute the leading truncation error terms component-wise using the equation\eqref{eqn:tau_linear} and the equation\eqref{eqn:tau_linear_RHS}. Once the computation was done, we made the comparison by evaluating $|\tau^{\text{RHS}}| - |\tau|$ at each near-boundary node and plotted the resulting differences as a contour over the $\theta_x$ and $\theta_y$ axes. Correspondingly, the numerical error differences $|\xi^{RHS}|-|\xi|$ between the RHS-affected and the original linear scheme case at the near-boundary nodes were also computed and compared. 

\begin{figure}[htbp]
	\centering
	\begin{subfigure}[b]{0.45\textwidth}
		\includegraphics[width=\linewidth]{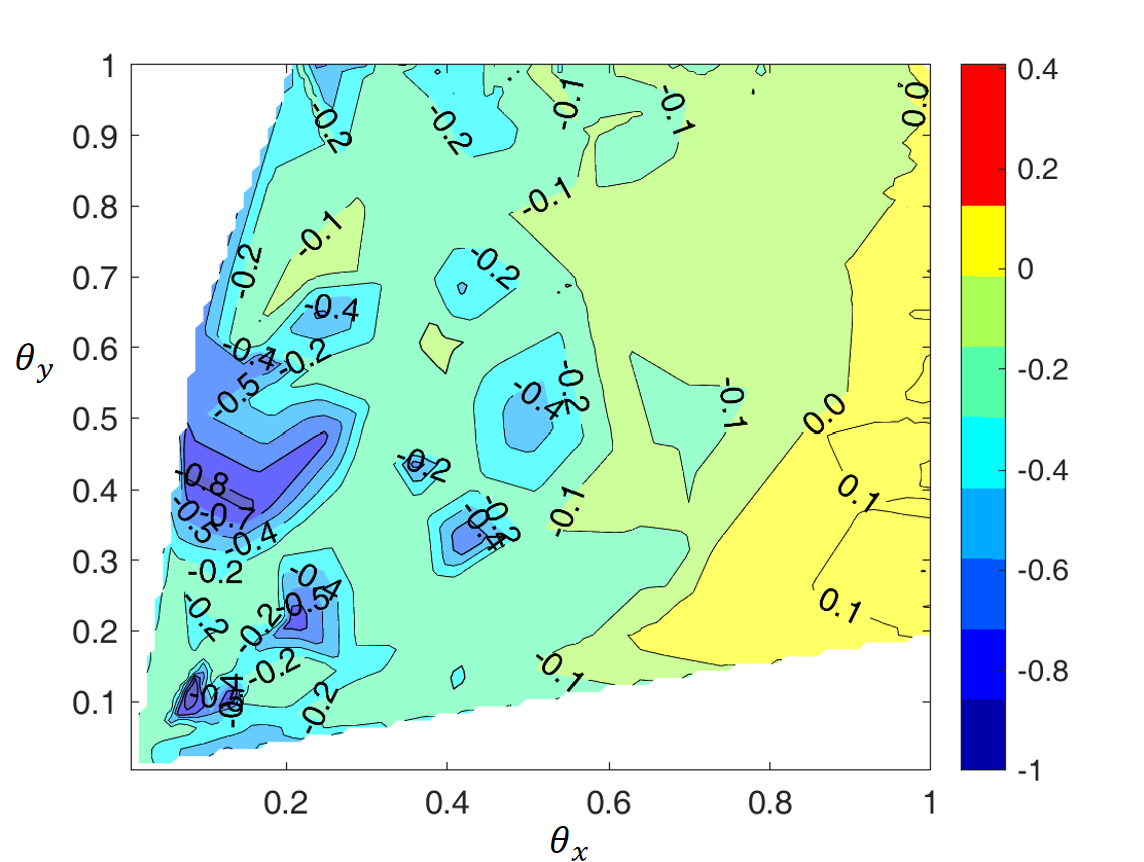}
		\caption{$\partial^2/\partial x^2$ part of $\left|\tau^{RHS}\right|$-$\left|\tau\right|$}
		\label{fig:linear_vs_linearRHS_2D_tauXpart}
	\end{subfigure}
	\hfill 
	\begin{subfigure}[b]{0.45\textwidth}
		\includegraphics[width=\linewidth]{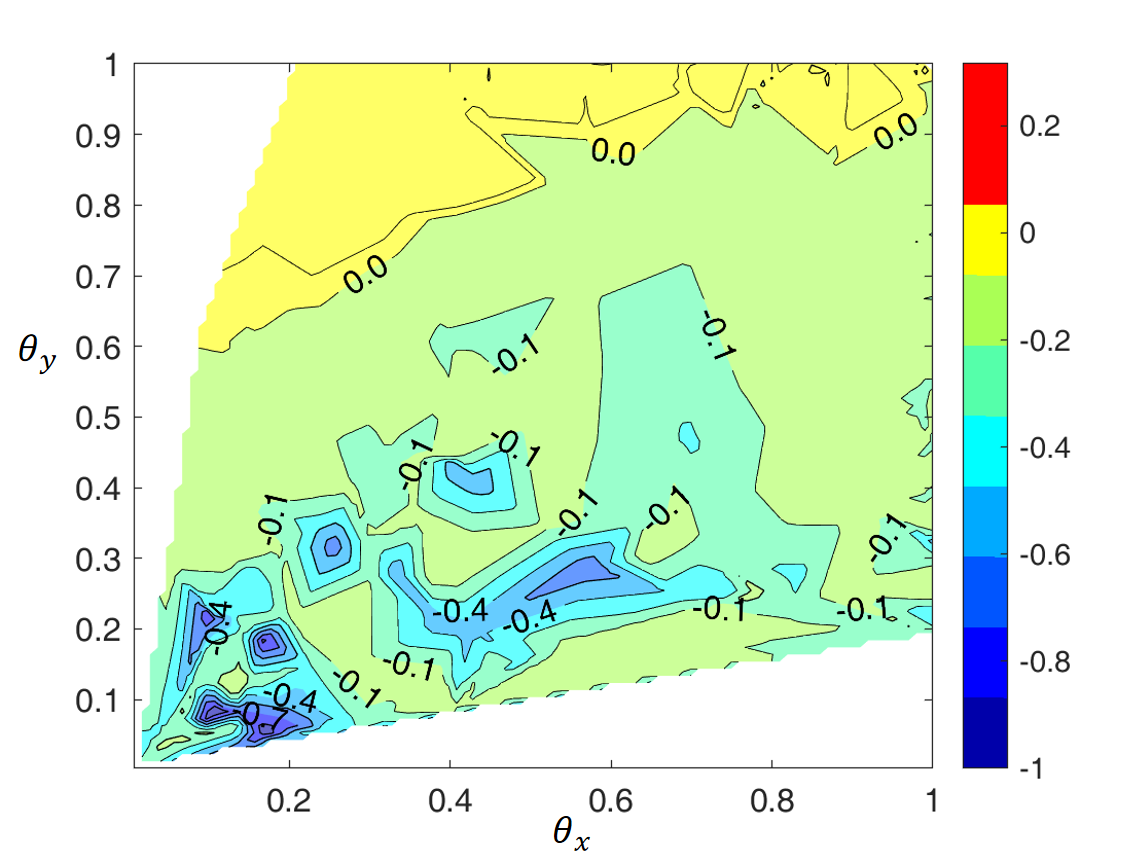}
		\caption{$\partial^2/\partial y^2$ part of $\left|\tau^{RHS}\right|$-$\left|\tau\right|$}
		\label{fig:linear_vs_linearRHS_2DtauYpart}
	\end{subfigure}
	
	\vspace{0.5cm} 
	
	\begin{subfigure}[b]{0.45\textwidth}
		\includegraphics[width=\linewidth]{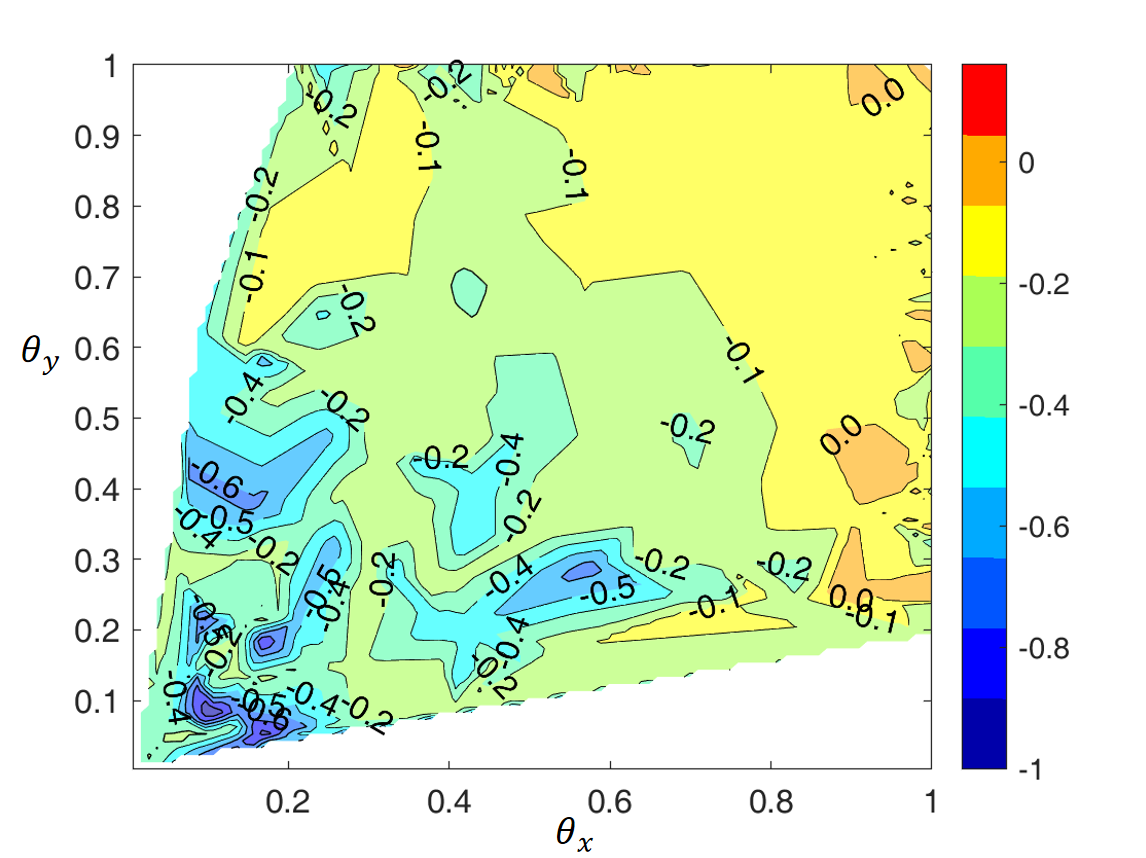}
		\caption{$\left|\tau^{RHS}\right|$-$\left|\tau\right|$}
		\label{fig:linear_vs_linearRHS_2DtauSum}
	\end{subfigure}
	\hfill
	\begin{subfigure}[b]{0.45\textwidth}
		\includegraphics[width=\linewidth]{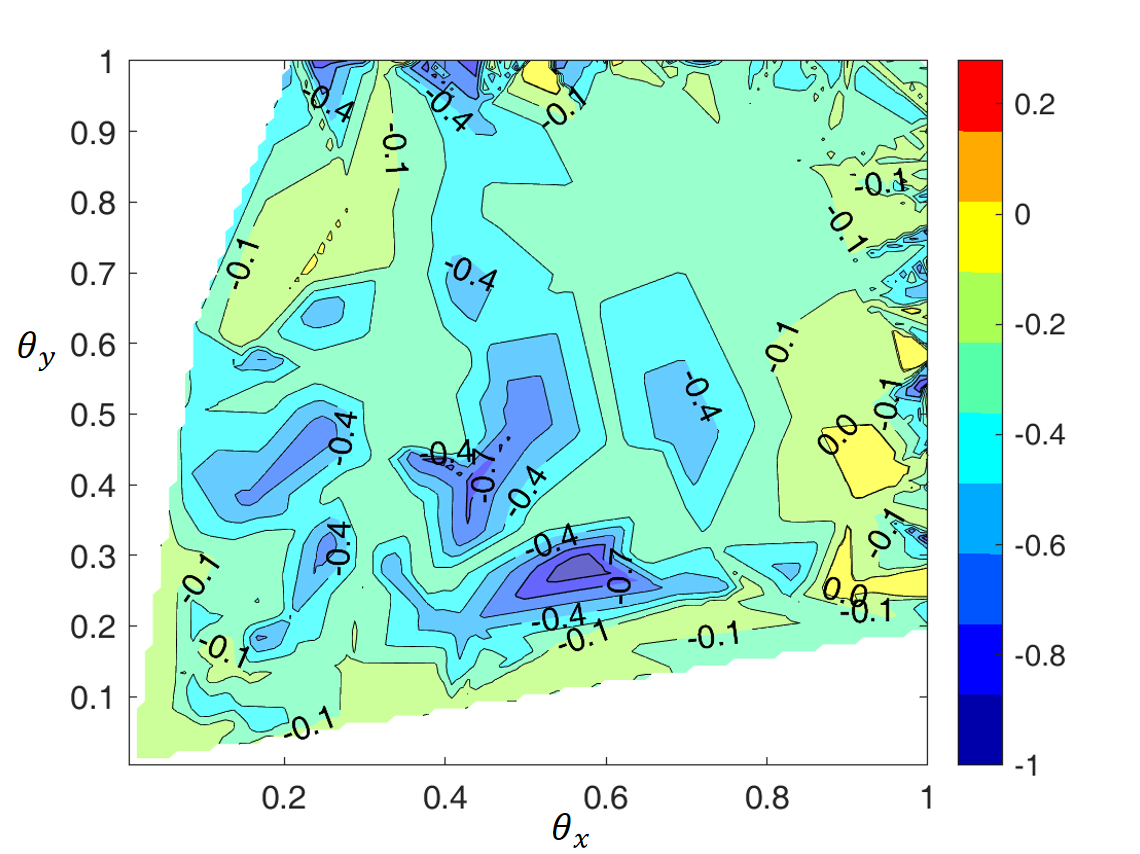}
		\caption{$\left|\xi^{RHS}\right|$-$\left|\xi\right|$}
		\label{fig:linear_vs_linearRHS_2Dxi}
	\end{subfigure}
	
	\caption{Contours of the comparison between the RHS-affected linear case and the ordinary linear case in terms of the truncation and numerical errors for all near-boundary nodes, each of which corresponds to a particular $(\theta_x,\theta_y)$ pair. A total number of $151\times151$ nodes were employed on the computational domain. (a) and (b) show the comparison of the leading truncation error component associated with $\partial\phi^2/\partial x^2$ and $\partial\phi^2/\partial y^2$, respectively. (c) compares the overall truncation error magnitude, and (d) compares the resulting numerical errors.}
	\label{fig:linear_vs_linearRHS_2D}
\end{figure}

Figure~\ref{fig:linear_vs_linearRHS_2D} presents these contour plots, where the results were normalized by their respective maximum absolute values. As shown in Figure~\ref{fig:linear_vs_linearRHS_2D_tauXpart}, for most $(\theta_x,\theta_y)$ pairs---except in the lower-right region of the $\theta$ domain---the RHS effect reduces the corresponding truncation error components (the $\partial^2/\partial x^2$ part). Comparing this contour with Figure~\ref{fig:linear_2D_tau_a} we observe very similar dividing lines (i.e., the zero level curves), which qualitatively validates the 2-D error analysis . Likewise, Figures~\ref{fig:linear_vs_linearRHS_2DtauYpart} and~\ref{fig:linear_vs_linearRHS_2DtauSum} exhibit patterns similar to those in Figures~\ref{fig:linear_2D_tau_b} and~\ref{fig:linear_2D_tau_c}, respectively.  Finally, Figure~\ref{fig:linear_vs_linearRHS_2Dxi} compares the numerical errors between the RHS-affected and original linear cases. Clearly, in most of the $\theta$ domain---except for small regions enclosed by the zero level curves in Figure~\ref{fig:linear_vs_linearRHS_2DtauSum}---the truncation errors for the RHS-affected linear scheme are smaller. Consequently, the resulting numerical errors also tend to be smaller, as seen from the similar zero-level region in Figure~\ref{fig:linear_vs_linearRHS_2Dxi}. These results suggest that the RHS effect generally improves the linear scheme by reducing the truncation errors at near-boundary nodes, and thus also reducing the numerical errors at those near-boundary nodes.

\begin{figure}[htbp]
	\centering
	\begin{subfigure}[b]{0.45\textwidth}
		\includegraphics[width=\linewidth]{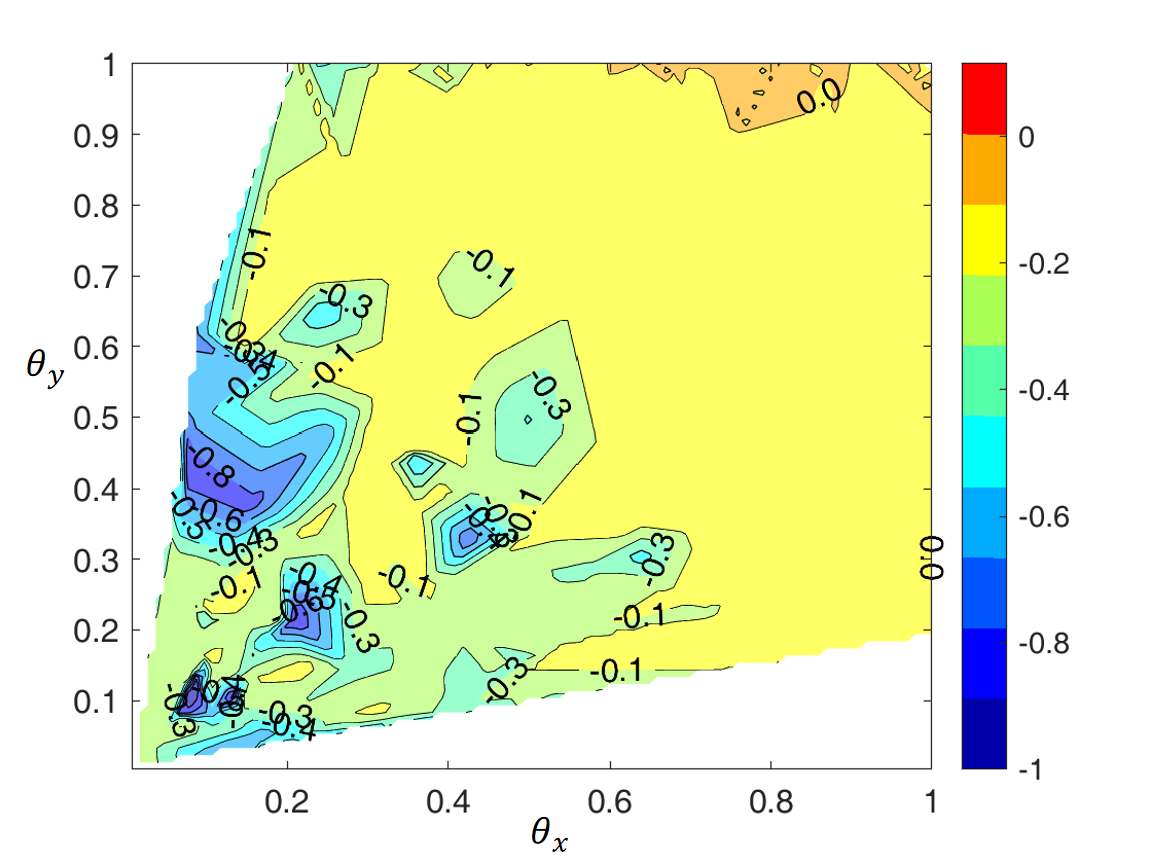}
		\caption{$\partial^2/\partial x^2$ part of $\left|\tau_{linear}\right|$-$\left|\tau_{quad}\right|$}
		\label{fig:linearRHS_vs_quadRHS_2DtauXpart}
	\end{subfigure}
	\hfill 
	\begin{subfigure}[b]{0.45\textwidth}
		\includegraphics[width=\linewidth]{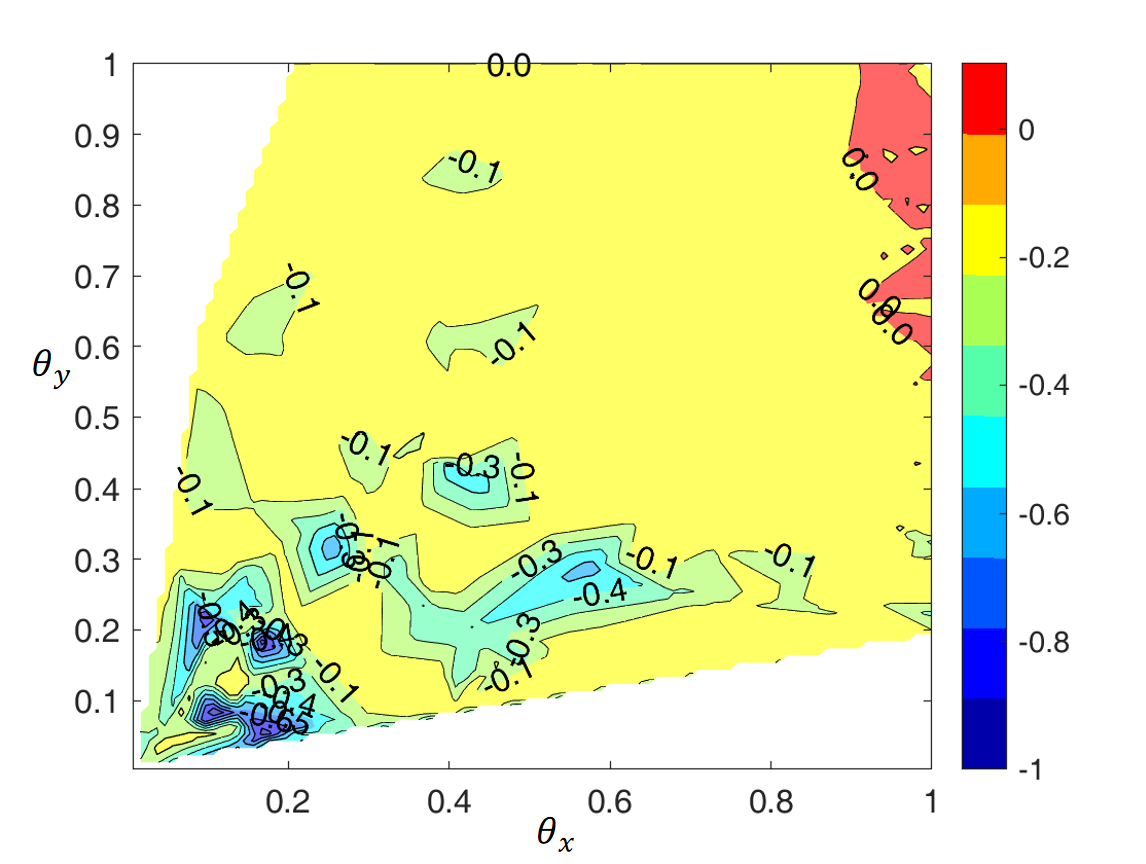}
		\caption{$\partial^2/\partial y^2$ part of $\left|\tau_{linear}\right|$-$\left|\tau_{quad}\right|$}
		\label{fig:linearRHS_vs_quadRHS_2DtauYpart}
	\end{subfigure}
	
	\vspace{0.5cm} 
	
	\begin{subfigure}[b]{0.45\textwidth}
		\includegraphics[width=\linewidth]{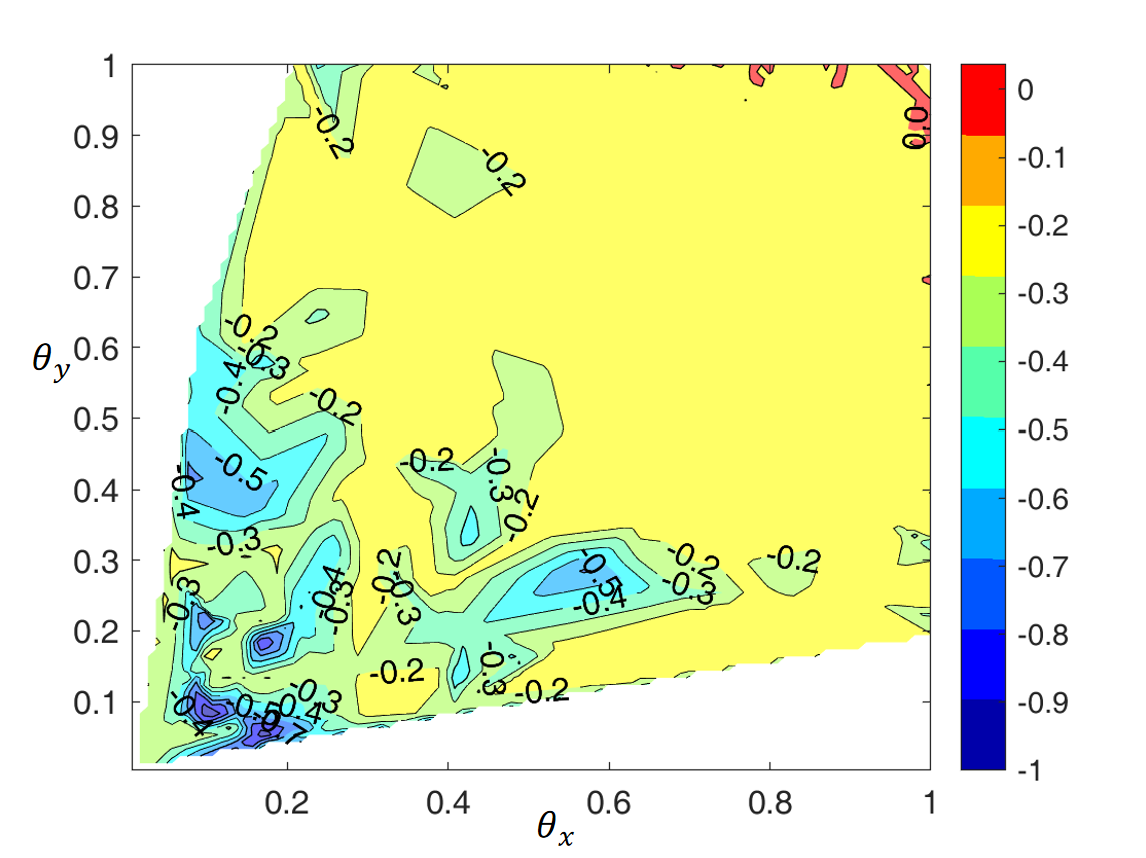}
		\caption{$\left|\tau_{linear}\right|$-$\left|\tau_{quad}\right|$}
		\label{fig:linearRHS_vs_quadRHS_2DtauSum}
	\end{subfigure}
	\hfill
	\begin{subfigure}[b]{0.45\textwidth}
		\includegraphics[width=\linewidth]{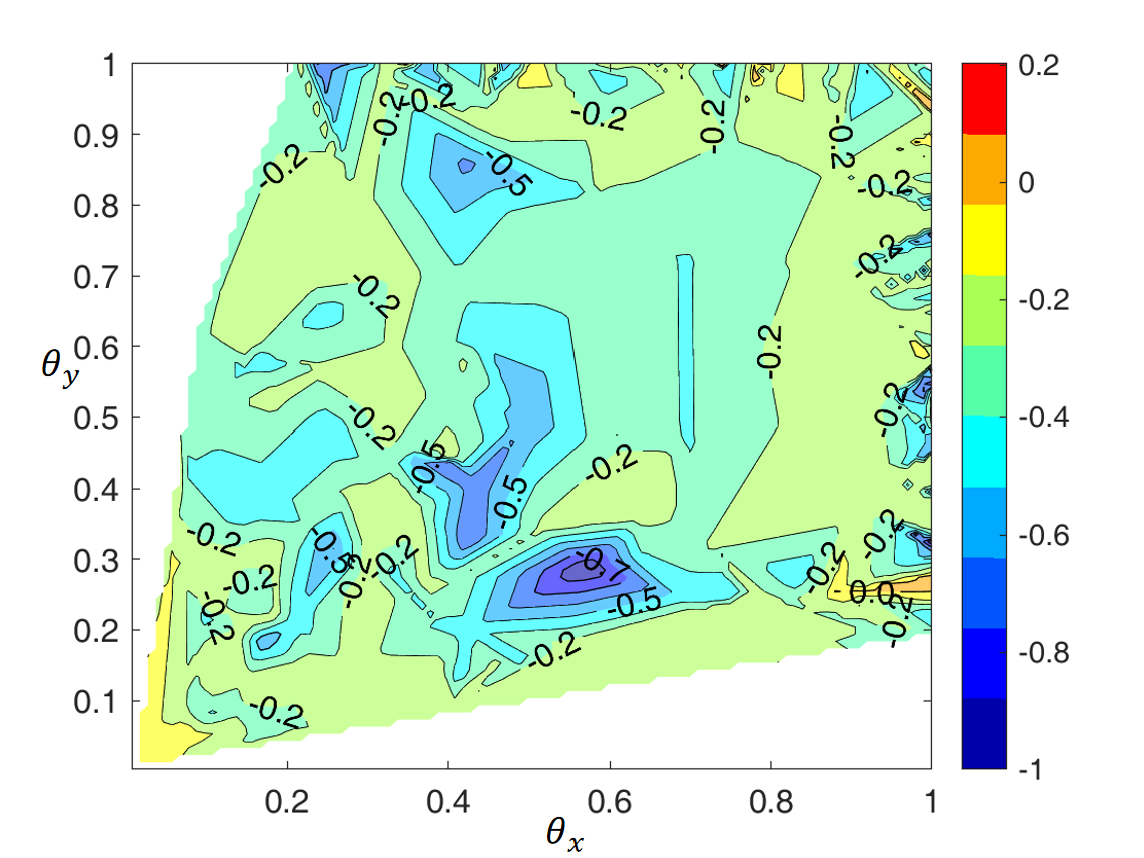}
		\caption{$\left|\xi_{linear}\right|$-$\left|\xi_{quad}\right|$}
		\label{fig:linearRHS_vs_quadRHS_2Dxi}
	\end{subfigure}
	
	\caption{Contours of the comparison between the RHS-affected linear case and the RHS-affected quad case in terms of the truncation and numerical errors for all near-boundary nodes, each of which corresponds to a particular $\theta_x$ and $\theta_y$ pair. (a) and (b) show the comparison of the leading truncation error component associated with $\partial\phi^2/\partial x^2$ and $\partial\phi^2/\partial y^2$, respectively. (c) compares the overall truncation error magnitude, and (d) compares the resulting numerical errors.}
	\label{fig:linearRHS_vs_quadRHS_2D}
\end{figure}

Next, we performed similar comparisons between the RHS-affected linear and RHS-affected quad cases. Figure~\ref{fig:linearRHS_vs_quadRHS_2D} presents these results. As seen from Figure~\ref{fig:linearRHS_vs_quadRHS_2DtauXpart}, \ref{fig:linearRHS_vs_quadRHS_2DtauYpart} and \ref{fig:linearRHS_vs_quadRHS_2DtauSum}, the RHS-affected linear scheme yields smaller truncation errors for almost all the $(\theta_x,\theta_y)$ pairs compared to the RHS-affected quadratic scheme. These findings are consistent with the earlier analysis presented in Figures~\ref{fig:LvsQ_2D_tau_a}, \ref{fig:LvsQ_2D_tau_b}, and~\ref{fig:LvsQ_2D_tau_c}, which showed that the linear scheme outperforms the quadratic scheme in terms of truncation errors at near-boundary nodes under RHS effects, assuming locally uniform RHS distribution and ``average" boundary shapes. It is also important to compare the resulting numerical errors at all near-boundary nodes, as shown in Figure~\ref{fig:linearRHS_vs_quadRHS_2Dxi}. Once again, the results support the conclusion that smaller truncation errors generally lead to smaller numerical errors at near-boundary nodes.

\subsubsection{Global Accuracy Assessment and RHS Error Compensation}
\label{subsec:overall_accuracy}
To quantify the numerical errors for each case, we measured the $L_1$ and $L_{\infty}$ norms of the numerical error vector $\Vec{\xi}$. The $L_1$ norms were normalized by the number of inner nodes to represent the “average” error of the numerical solution, whereas the $L_{\infty}$ norms represent the maximum errors. Using error results obtained from using different node spacing $\Delta_1$ and $\Delta_2$, we could estimate the error order $p$, using the formula $p=\ln(L_{\Delta_1}/L_{\Delta_2})/\ln(\Delta_2/\Delta_1)$. For each case, we listed these quantities in table form, and then juxtaposed the plots of the truncation errors and the resulting numerical errors at all nodes.

We first demonstrate the accuracy results from solving the  2-D problem by both the linear and the quadratic schemes, without the inaccurate RHS effect. Table \ref{tab:2D accuracy} summarizes the numerical error results, tested on different meshes. It is seen that when the number of nodes were increased, both the linear and quadratic scheme demonstrated second-order accuracy, although the errors were consistently lower for the quadratic setup, especially when measured with $L_{\infty}$ norms. However, when the number of nodes reached $161\times161$, the solution quality suddenly dropped for the quadratic case, despite the finer mesh grids. Further investigation revealed that this phenomenon occurs whenever an inner node is situated between two Dirichlet boundary interfaces along a specific axial direction (i.e., it has no neighboring nodes along that direction). In such cases, while the quadratic scheme remains applicable, it generates substantial truncation errors at these nodes, thereby degrading the overall solution accuracy. To avoid this issue, some detection measure must be employed, and the treatment for these "troubled" nodes should fall back to the linear scheme to avoid failure due to insufficient stencil points. Once corrected, the resulting $L_1$ and $L_{\infty}$ norms were reduced to $1.508\times10^{-8}$ and $6.131\times10^{-8}$, respectively. This indicates that the quadratic scheme is less robust than the linear one. For this reason, later (and earlier) numerical error plots were generated on a $151\times151$ mesh, where such situation is absent.

\begin{table}[]
	\caption{numerical errors from solving the 2-D case, using accurate RHS values $b_{i,j}$}
	\label{tab:2D accuracy}
	\resizebox{\textwidth}{!}{%
		\begin{tabular}{l|llll|llll}
			\hline
			& \multicolumn{4}{c|}{linear extrapolation}        & \multicolumn{4}{c}{quadratic extrapolation}            \\ \hline
			\begin{tabular}[c]{@{}l@{}}number \\ of nodes\end{tabular} &
			$L_1$ error &
			order &
			$L_\infty$ error&
			order &
			$L_1$ error &
			order &
			$L_\infty$ error&
			order \\
			41$\times$41 &
			8.013e-6 &
			\multicolumn{1}{c}{--} &
			2.047e-5 &
			\multicolumn{1}{c|}{--} &
			2.3643e-7 &
			\multicolumn{1}{c}{--} &
			6.7915e-7 &
			\multicolumn{1}{c}{--} \\
			81$\times$81   & 2.0971e-6 & 1.93 & 6.198e-6 & 1.72 & 5.8779e-8 & 2.00 & 1.7224e-7 & 1.82 \\
			161$\times$161 & 5.262e-7  & 1.99 & 1.654e-6 & 1.91 & \textbf{1.1756e-6} & \textbf{-4.3}  & \textbf{9.7126e-4} & \textbf{-12.5}  \\ \hline
		\end{tabular}%
	}
\end{table}

\begin{figure}
	\begin{subfigure}{0.49\textwidth}
		\includegraphics[width=\linewidth]{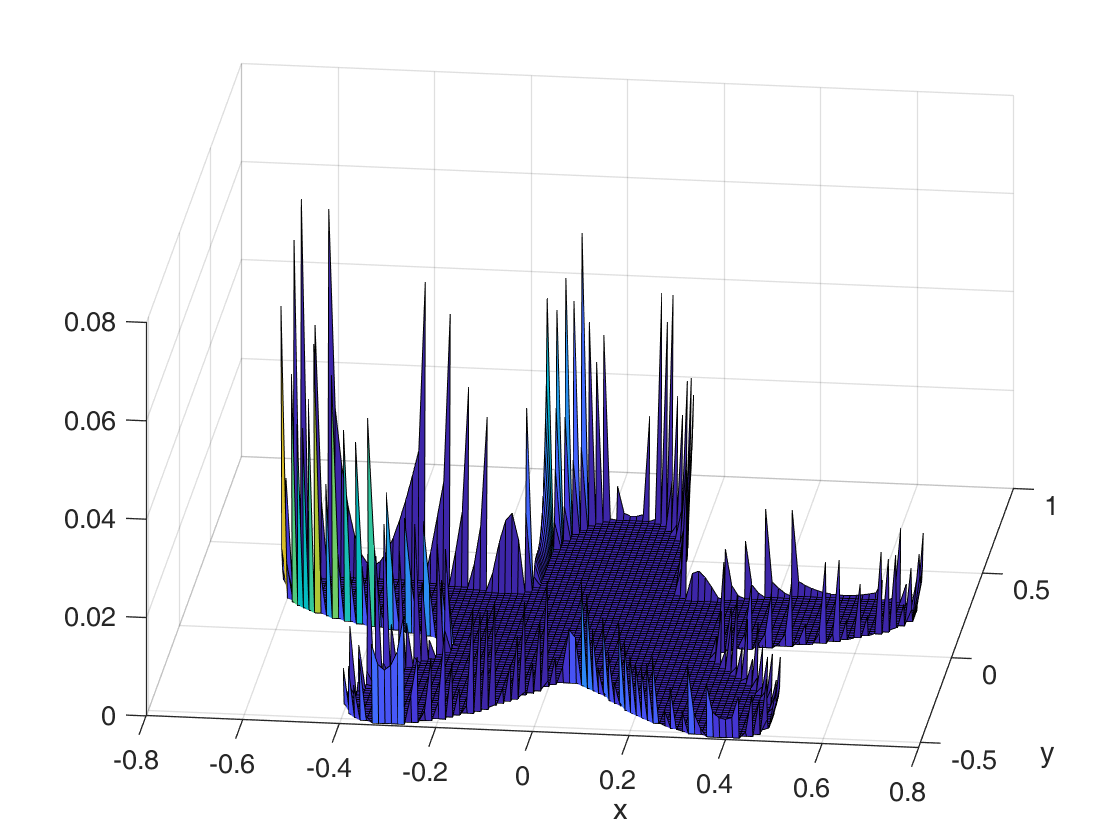}
		\caption{truncation errors $\tau$}
		\label{fig:linear_surf_a}
	\end{subfigure}%
	\hspace*{\fill}   
	\begin{subfigure}{0.49\textwidth}
		\includegraphics[width=\linewidth]{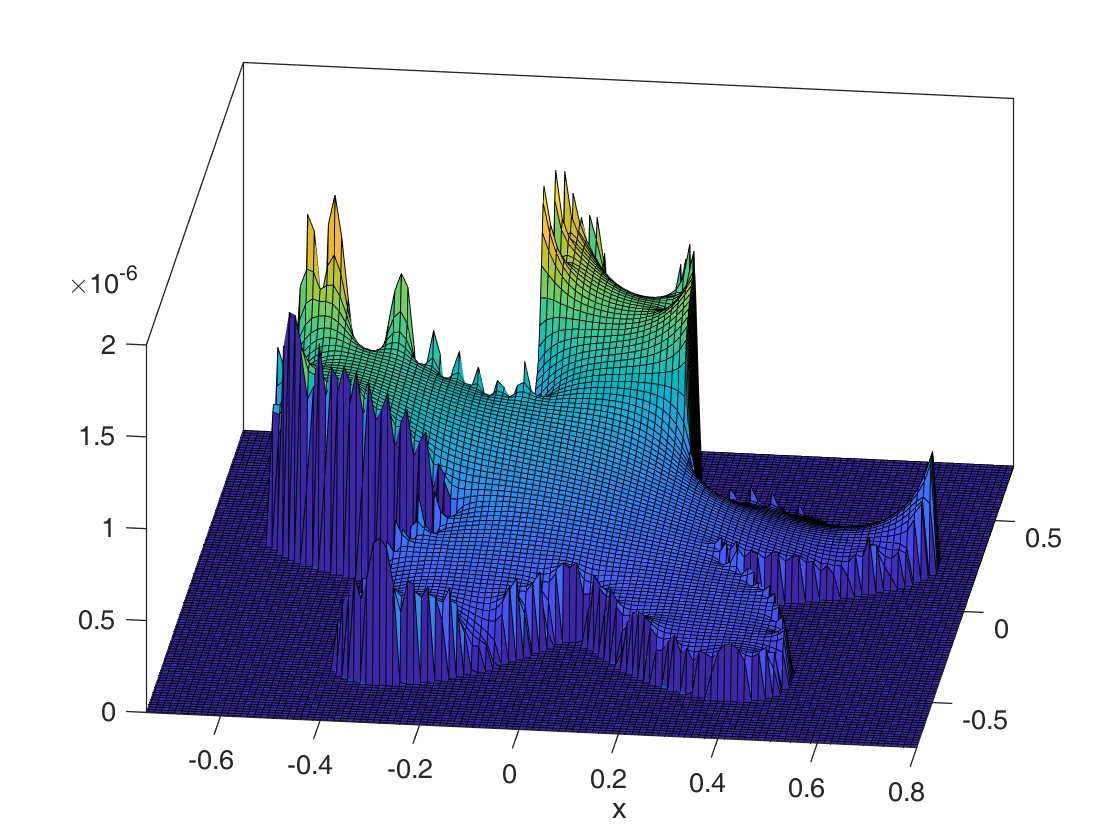}
		\caption{numerical errors $\xi$ (inverted)} 
		\label{fig:linear_surf_b}
	\end{subfigure}
	
	\caption{The plots of the truncation errors and the final numerical errors at all nodes from solving the 2-D Poisson equation by the linear scheme. (a) shows the truncation errors, and (b) shows the resulting numerical errors, which were inverted for better comparison.} 
	\label{fig:linear_surf}
\end{figure}

\begin{figure}
	\begin{subfigure}{0.49\textwidth}
		\includegraphics[width=\linewidth]{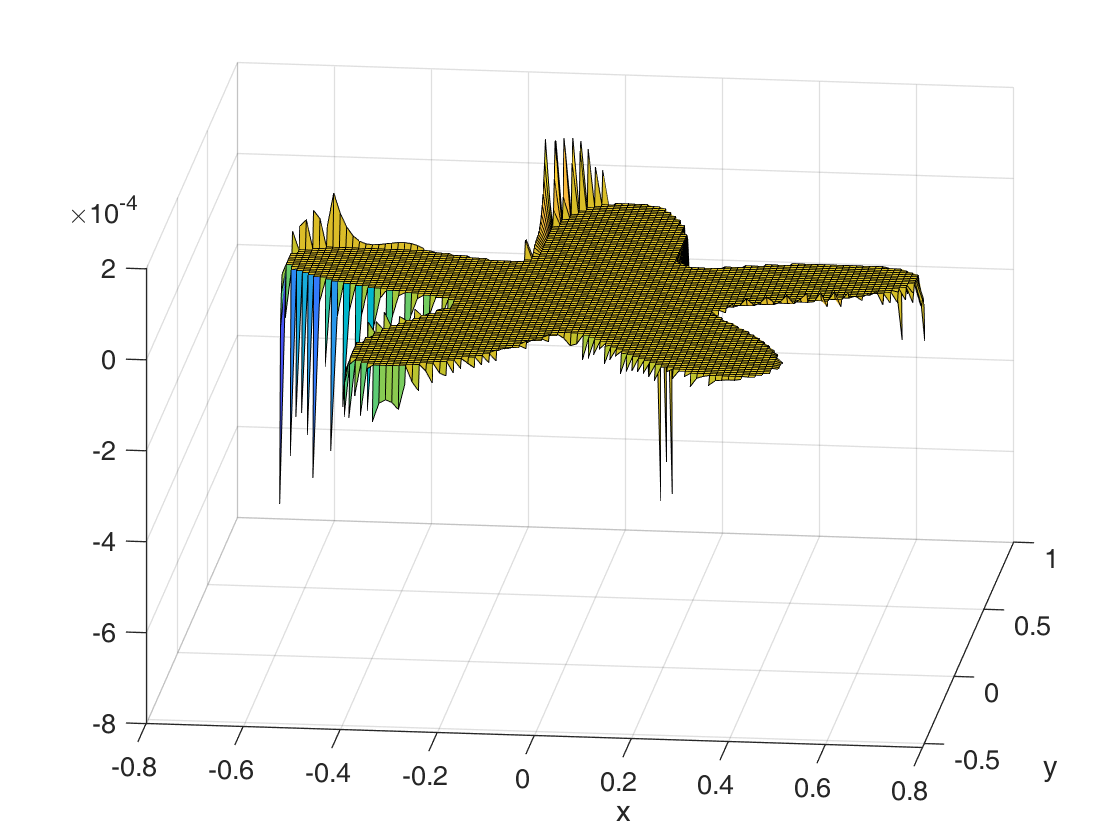}
		\caption{truncation errors $\tau$}
		\label{fig:quad_surf_a}
	\end{subfigure}%
	\hspace*{\fill}   
	\begin{subfigure}{0.49\textwidth}
		\includegraphics[width=\linewidth]{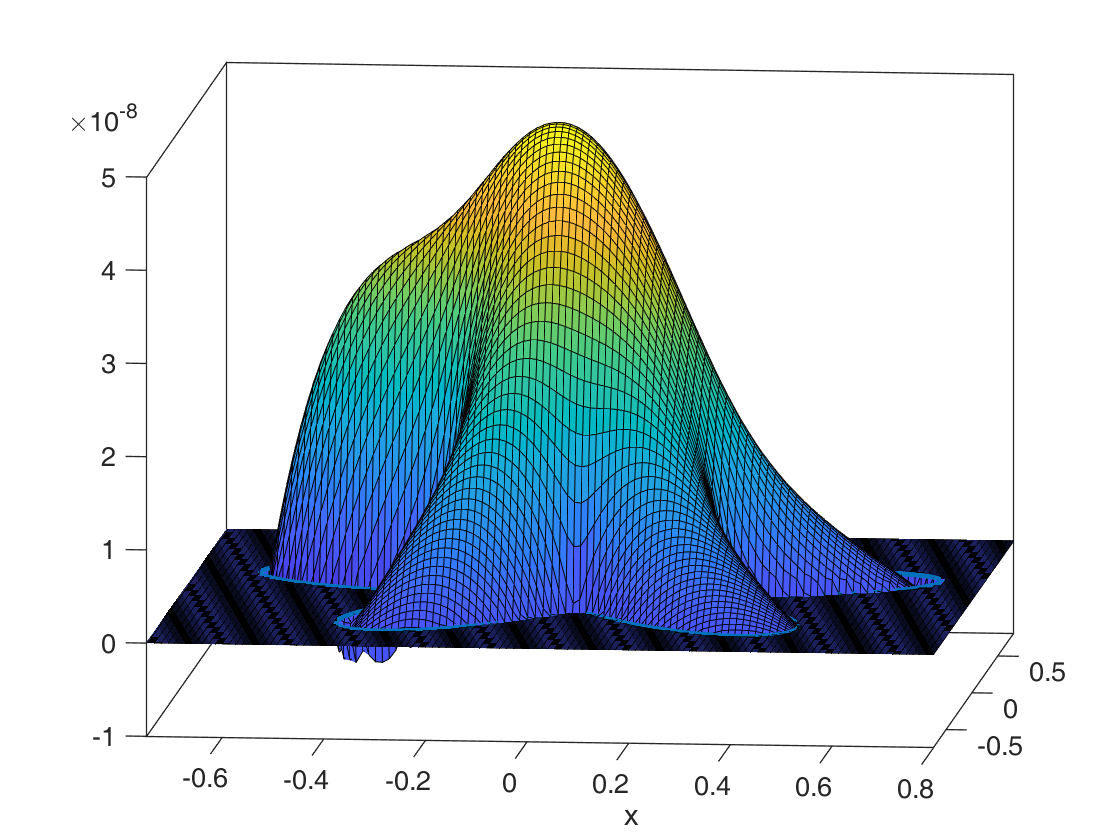}
		\caption{numerical errors $\xi$ (inverted)} 
		\label{fig:quad_surf_b}
	\end{subfigure}
	
	\caption{The plots of the truncation errors and the final numerical errors at all nodes from solving the 2-D Poisson equation by the quadratic scheme. (a) shows the truncation errors, and (b) shows the resulting numerical errors, which were inverted for better comparison.} 
	\label{fig:quad_surf}
\end{figure}

Figure~\ref{fig:linear_surf} presents the truncation and numerical errors at all nodes obtained from solving the 2-D Poisson equation using the linear scheme with accurate RHS values. As shown in Figure~\ref{fig:linear_surf_a}, the truncation errors near the boundary—being zeroth-order—are significantly larger than those in the interior, which are second-order. Consequently, the numerical errors in Figure~\ref{fig:linear_surf_b} are most pronounced at the boundaries and decrease smoothly in magnitude toward the interior of the domain. Figure~\ref{fig:quad_surf} shows the corresponding errors obtained using the quadratic scheme with accurate RHS values. As illustrated in figure~\ref{fig:quad_surf_a}, although the truncation errors--now first-order--remain larger near the boundaries, their magnitudes are substantially reduced compared to those in the linear scheme. This reduction results in a qualitatively different numerical error distribution, which now exhibits a center-dominant pattern, as seen in figure~\ref{fig:quad_surf_b}. It is important to note that, based on previous 1-D error analysis showing that boundary truncation errors and numerical errors have opposite signs, the numerical error plots have been inverted here to facilitate direct visual comparison with the corresponding truncation error plots.                                                                                                                                                                               

Next, we present the results for the two-dimensional case in the presence of the inaccurate right-hand side (RHS) effect. The numerical errors are summarized in Table \ref{tab:2D accuracy_RHS}. It can be observed that both the linear and quadratic schemes retain second-order accuracy. However, upon examining the $L_1$ and $L_\infty$ error magnitudes between the two schemes, it could seen that the RHS-affected linear scheme consistently outperformed those of the quadratic counterpart. Moreover, when these results are compared with those in Table\ref{tab:2D accuracy}, it becomes evident that the linear scheme benefits from the use of the inaccurate RHS values $\overline{b_{i,j}}$, whereas the quadratic scheme experiences a notable degradation in accuracy. These findings extend the previous analysis by demonstrating that the use of $\overline{b_{i,j}}$ not only alters the boundary truncation errors $\tau$, and consequently the boundary numerical errors $\xi$, but impacts the overall accuracy across the entire computational domain.

\begin{figure}
	\begin{subfigure}{0.49\textwidth}
		\includegraphics[width=\linewidth]{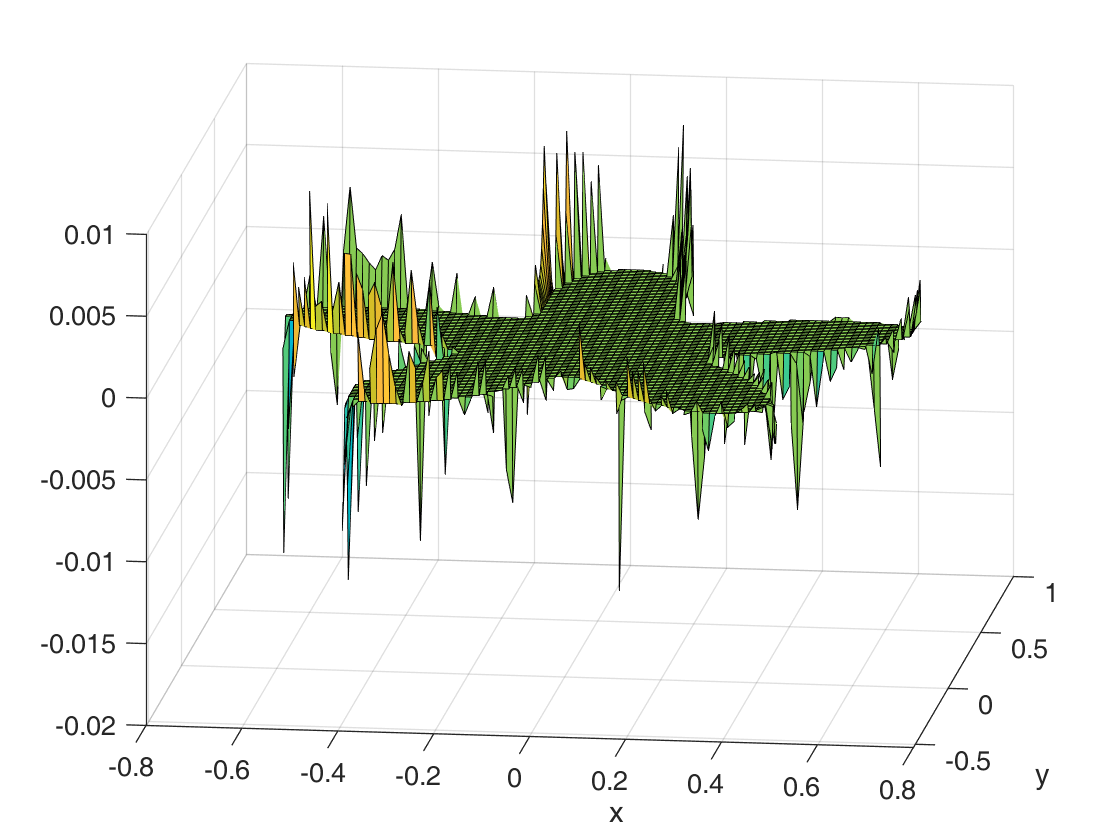}
		\caption{truncation errors $\tau$}
		\label{fig:linear_surf_RHS_a}
	\end{subfigure}%
	\hspace*{\fill}   
	\begin{subfigure}{0.49\textwidth}
		\includegraphics[width=\linewidth]{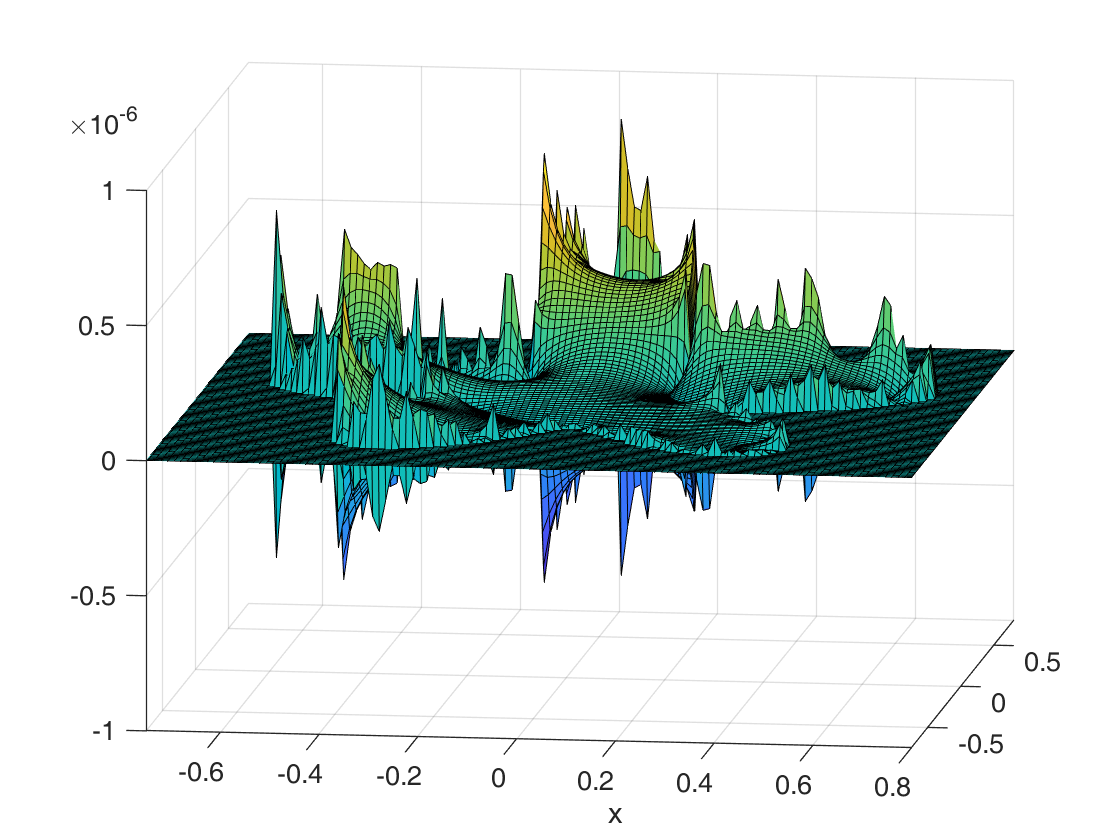}
		\caption{numerical errors $\xi$ (inverted)}
		\label{fig:linear_surf_RHS_b} 
	\end{subfigure}
	
	\caption{The plot of the truncation errors and the final numerical errors at all nodes from solving the 2-D Poisson equation by the linear scheme with the inaccurate RHS effect. (a) shows the truncation errors, and (b) shows the resulting numerical errors, which were inverted for better comparison.} 
	\label{fig:linear_surf_RHS}
\end{figure}

Figure~\ref{fig:linear_surf_RHS} presents the truncation and numerical errors associated with solving the 2-D Poisson equation using the linear scheme with inaccurate right-hand side (RHS) values $\overline{b_{i,j}}$. As illustrated in Figure~\ref{fig:linear_surf_RHS_a}, the truncation errors near the boundaries remain dominant relative to those in the interior; however, their magnitudes are reduced compared to the case with accurate RHS values (see Figure~\ref{fig:linear_surf_a}). Notably, the signs of the truncation errors near the boundaries may change at certain nodes. Correspondingly, the numerical errors, shown in Figure~\ref{fig:linear_surf_RHS_b}, also exhibit an overall reduction in magnitude and sign changes in specific regions.

\begin{figure}
	\begin{subfigure}{0.49\textwidth}
		\includegraphics[width=\linewidth]{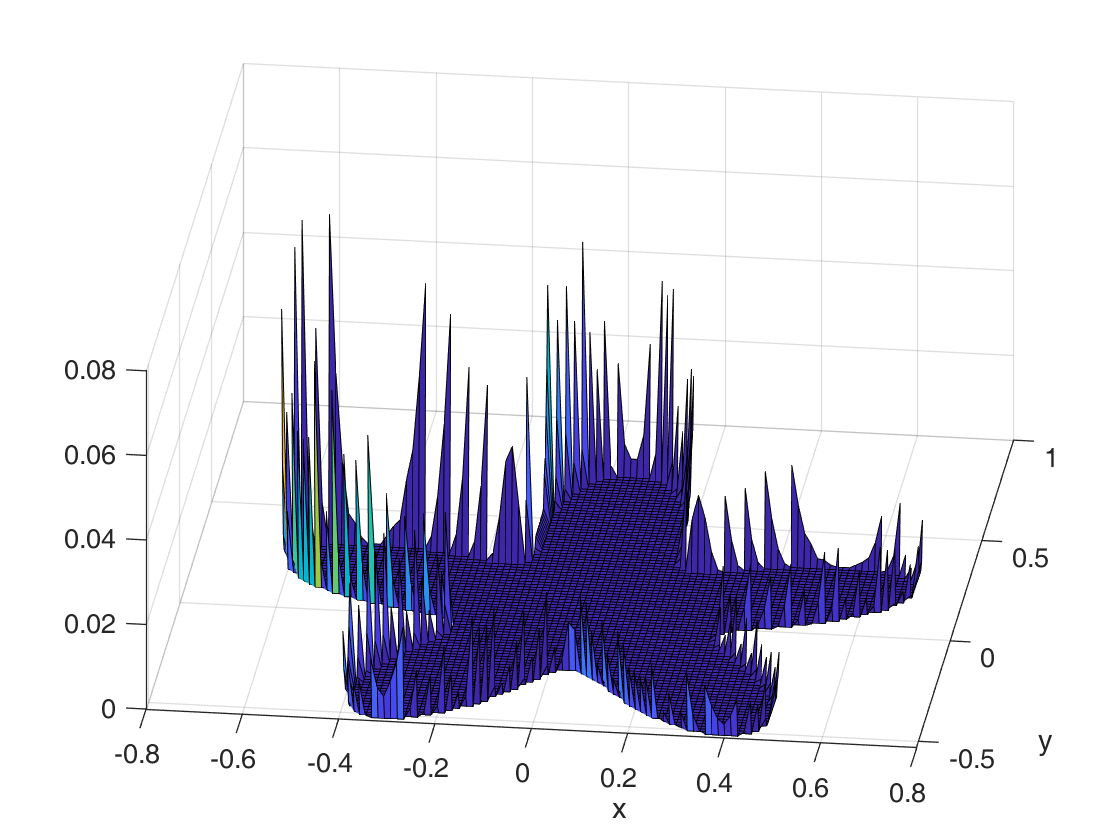}
		\caption{truncation errors $\tau$ (inverted)}
		\label{fig:quad_surf_RHS_a}
	\end{subfigure}%
	\hspace*{\fill}   
	\begin{subfigure}{0.49\textwidth}
		\includegraphics[width=\linewidth]{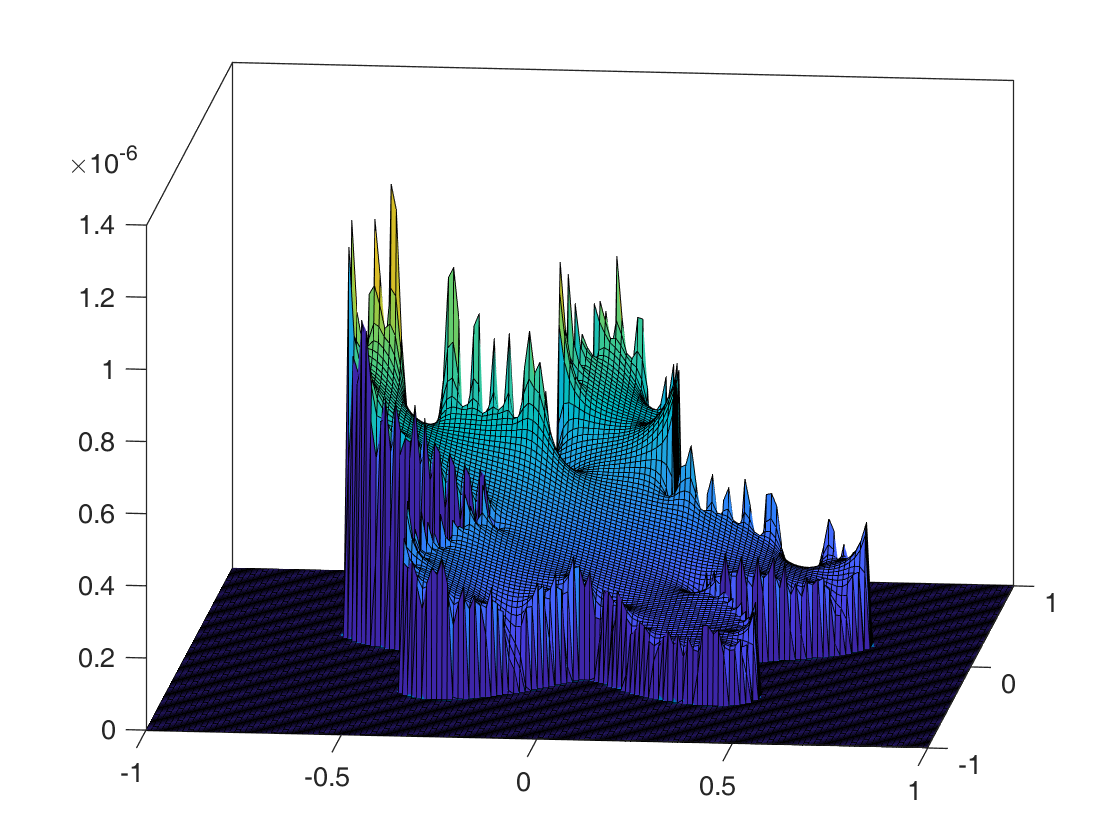}
		\caption{numerical errors $\xi$}
		\label{fig:quad_surf_RHS_b} 
	\end{subfigure}
	
	\caption{The plot of the truncation errors and the final numerical errors at all nodes from solving the 2-D Poisson equation by the quadratic scheme with the inaccurate RHS effect. (a) shows the truncation errors, which were inverted for better comparison, and (b) shows the resulting numerical errors.} 
	\label{fig:quad_surf_RHS}
\end{figure}

Figure~\ref{fig:quad_surf_RHS} displays the corresponding error plots obtained using the quadratic scheme with the same inaccurate RHS values. For clarity, the truncation error plot ($\tau$) is inverted instead of the numerical error plot ($\xi$). As evident from Figure~\ref{fig:quad_surf_RHS_a}, the truncation errors near the boundary are approximately two orders of magnitude larger than those observed in Figure~\ref{fig:quad_surf_a}, where accurate RHS values were used. This increase aligns with the earlier 2-D error analysis, which indicates that inaccurate RHS values modify the boundary truncation error from a first-order term to a zeroth-order term. Given that the numerical experiments used a grid size of $2/150 \approx 0.01$, the observed difference in magnitude is consistent with theoretical expectations.
\begin{table}[]
	\caption{numerical errors from solving the 2-D case, using inaccurate RHS values $\overline{b_{i,j}}$}
	\label{tab:2D accuracy_RHS}
	\resizebox{\textwidth}{!}{%
		\begin{tabular}{l|llll|llll}
			\hline
			& \multicolumn{4}{c|}{linear extrapolation}        & \multicolumn{4}{c}{quadratic extrapolation}            \\ \hline
			\begin{tabular}[c]{@{}l@{}}number \\ of nodes\end{tabular} &
			$L_1$ error &
			order &
			$L_\infty$ error&
			order &
			$L_1$ error &
			order &
			$L_\infty$ error&
			order \\
			41$\times$41 &
			1.3931e-6 &
			\multicolumn{1}{c}{--} &
			7.7411e-6 &
			\multicolumn{1}{c|}{--} &
			5.1160e-6 &
			\multicolumn{1}{c}{--} &
			1.7039e-5 &
			\multicolumn{1}{c}{--} \\
			81$\times$81   & 3.8876e-7 & 1.84 & 2.8273e-6 & 1.45 & 1.3768e-6 & 1.89 & 4.2559e-6 & 2.00 \\
			151$\times$151 & 1.0415e-7  & 2.10 & 8.4515e-7 & 1.92 & 3.6962e-7 & 2.09  & 1.2271e-6 & 1.98  \\ \hline
		\end{tabular}%
	}
\end{table}

Finally, we tested the effectiveness of using the approximated $\overline{\delta_{i,j}}$ values to counteract the effect of inaccurate RHS values. It is important to note that the exact $\delta_{i,j}$ values---defined as the ratio between the inaccurate and accurate RHS values---can be computed as $\delta_{i,j} = \overline{b_{i,j}} / b_{i,j}$, since the exact RHS values $b_{i,j}$ were known in our numerical tests. In contrast, the approximated values $\overline{\delta_{i,j}}$ were derived under the assumption of a locally uniform distribution of $b(x,y)$ and were evaluated using equation\eqref{eqn:delta}. Consequently, we began our tests on the 2-D case while varying the exact solution to $\phi(x,y)=x^2+y^2$, which corresponds to a uniform RHS distribution $b(x,y)=4$. Three testing conditions were considered: in the first, accurate RHS values were used for all nodes; in the second, the RHS values for the near-boundary nodes were computed by equation \ref{eqn:rhs_calc}; in the third, the RHS values for the near-boundary nodes were modified by $\widetilde{b_{i,j}}=\overline{b_{i,j}}/\overline{\delta_{i,j}}$, in an effort to calibrate them to closely approximate the true $b_{i,j}$. The resulting numerical errors under each condition are presented in Table~\ref{tab:2D rhs effect uniform}, where the $L_{\infty}$ norm was used for evaluation. In this special case, the quadratic scheme demonstrates pronounced advantage, reaching an numerical error in the order of $10^{-16}$ while using only $41\times41$ nodes, provided that accurate RHS value were used. This extreme level of accuracy is due to the fact that the exact solution of $\phi(x,y)$ was given in quadratic form, in which case the quadratic scheme provides the perfect discretization of the LHS of the Poisson equation. Interestingly, further refining the mesh to $81\times81$ resulted in a slight degradation of accuracy due to the accumulation of rounding errors during the solution of the linear system. When $161\times161$ nodes were used, the solution from the quadratic scheme deteriorated sharply, obviously due to the rare situation that certain nodes were located between Dirichlet boundary interfaces, as discussed previously. Under the second test condition where inaccurate RHS values $\overline{b_{i,j}}$ were used for the near-boundary nodes, the accuracy of the linear scheme improved, whereas that of the quadratic scheme decline---consistent with earlier observations. Finally, in the third condition, the numerical results perfectly matched those obtained using accurate RHS values. This was expected, as the approximated $\overline{\delta_{i,j}}$ values were highly accurate in this special case, given the true uniformity of the RHS distribution $b(x,y)=4$.
\begin{table}[]
	\caption{numerical errors from solving the varied 2-D case, where the RHS values are uniformly distributed as $b=4$.}
	\label{tab:2D rhs effect uniform}
	\setlength{\extrarowheight}{5pt} 
	\resizebox{\textwidth}{!}{%
		\begin{tabular}{l|cc|cc|cc}
			\hline
			& \multicolumn{2}{c|}{Accurate RHS $b$} & \multicolumn{2}{c|}{Inaccurate RHS $\overline{b_{i,j}}$} & \multicolumn{2}{c}{Modified $\widetilde{b_{i,j}}=\overline{b_{i,j}}/\overline{\delta_{i,j}}$} \\ \hline
			\begin{tabular}[c]{@{}l@{}}number of \\ nodes used\end{tabular} &
			\multicolumn{1}{l}{Linear $L_{\infty}$} &
			\multicolumn{1}{l|}{Quad $L_{\infty}$} &
			\multicolumn{1}{l}{{Linear $L_{\infty}$}} &
			\multicolumn{1}{l|}{{Quad $L_{\infty}$}} &
			\multicolumn{1}{l}{{Linear $L_{\infty}$}} &
			\multicolumn{1}{l}{{Quad $L_{\infty}$}} \\ \hline
			41$\times$41   & 5.728e-4      & 2.498e-16     & 3.499e-4    & 6.201e-4   & 5.708e-4     & 2.498e-16     \\
			81$\times$81   & 1.446e-4      & 5.551e-16     & 8.386e-5    & 1.482e-4   & 1.446e-4     & 5.551e-16     \\
			\textbf{161$\times$161} & 3.803e-5      & \textbf{1.032e-6}      & 1.652e-5    & \textbf{3.603e-5}   & 3.803e-5     & \textbf{1.032e-6}      \\ \hline
		\end{tabular}
	}
\end{table}

\begin{figure}
	\begin{subfigure}{0.49\textwidth}
		\includegraphics[width=\linewidth]{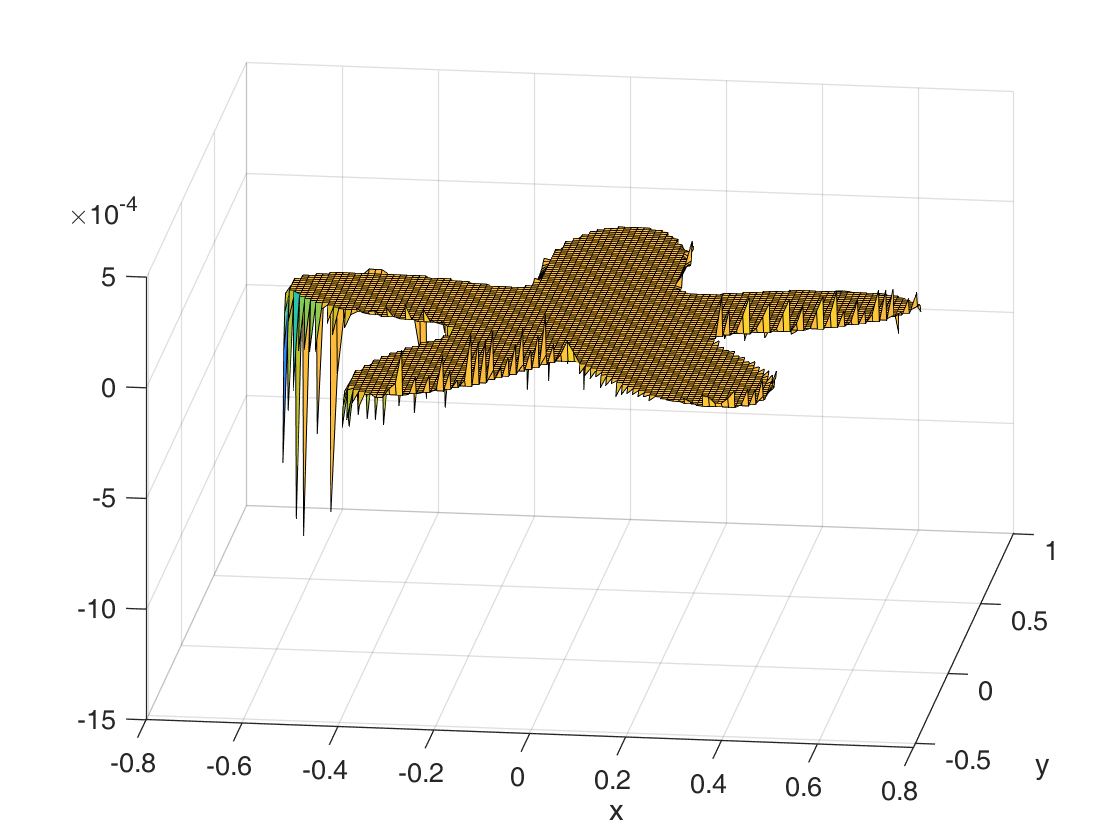}
		\caption{truncation errors $\tau$}
	\end{subfigure}%
	\hspace*{\fill}   
	\begin{subfigure}{0.49\textwidth}
		\includegraphics[width=\linewidth]{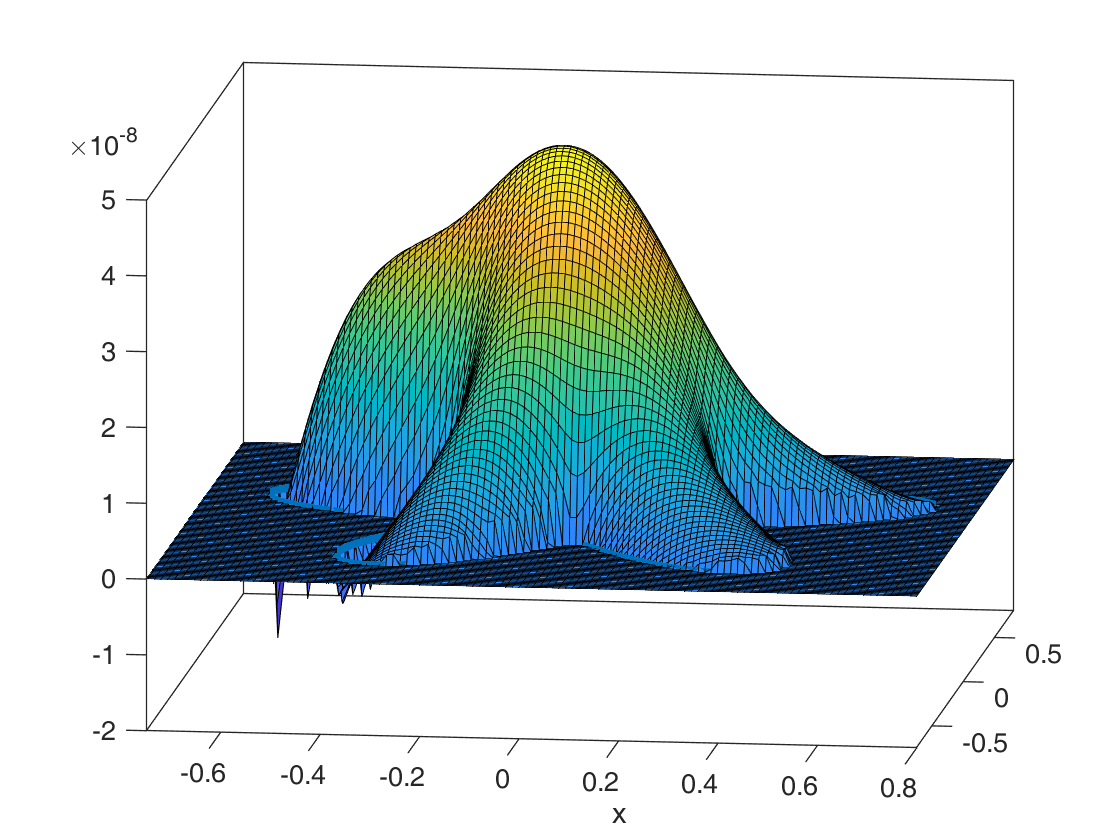}
		\caption{numerical errors $\xi$ (inverted)} 
	\end{subfigure}
	
	\caption{The plot of the truncation errors and the final numerical errors at all nodes from solving the 2-D Poisson equation by the RHS-affected quadratic scheme, where RHS values were modified by $b_{i,j}=\overline{b_{i,j}}/\overline{\delta_{i,j}}$ for the near-boundary nodes. (a) shows the truncation errors, and (b) shows the resulting numerical errors, which were inverted for better comparison.} 
	\label{fig:quad_surf_FIX}
\end{figure}

We then repeated the same procedures for the original 2-D case, where the exact solution was given by $\phi(x,y) = \left[(x+2)^2 + (y-2)^2\right]^{-1}$, corresponding to a nonuniform RHS distribution. The results are presented in Table~\ref{tab:2D rhs effect nonuniform}. As expected, under the first condition with accurate RHS values, the quadratic scheme outperformed the linear scheme, consistent with its higher-order accuracy of the truncation errors at near-boundary nodes. In the second condition, where inaccurate RHS values $\overline{b_{i,j}}$ were applied to near-boundary nodes, the solution from the linear scheme improved, whereas the performance of the quadratic scheme declined—an outcome that aligns with earlier observations. Crucially, in the third condition, the application of approximated coefficients $\overline{\delta_{i,j}}$ to calibrate the inaccurate RHS values (i.e., using $\widetilde{b_{i,j}} = \overline{b_{i,j}} / \overline{\delta_{i,j}}$ for near-boundary nodes) resulted in accuracy comparable to that obtained using the true RHS values. Interestingly, the quadratic scheme yielded even slightly lower errors than those from the second condition. This demonstrates that, despite $\overline{\delta_{i,j}}$ being derived under the assumption of local uniformity in $b(x,y)$, it closely approximates the true $\delta_{i,j}$ values even in nonuniform cases. Figure~\ref{fig:quad_surf_FIX} shows the corresponding error plots for the quadratic scheme using the calibrated RHS values, and demonstrates strong visual agreement with Figure~\ref{fig:quad_surf}, where accurate RHS values were applied across the entire domain.

\begin{table}[]
	\caption{2-D test case with nonuniform $b_{i,j}$---the numerical error results obtained using different treatments of RHS values.}
	\label{tab:2D rhs effect nonuniform}
	\setlength{\extrarowheight}{5pt} 
	\resizebox{\textwidth}{!}{%
		\begin{tabular}{l|cc|cc|cc}
			\hline
			& \multicolumn{2}{c|}{Accurate RHS $b_{i,j}$} & \multicolumn{2}{c|}{Inaccurate RHS $\overline{b_{i,j}}$} & \multicolumn{2}{c}{Modified $\widetilde{b_{i,j}}=\overline{b_{i,j}}/\overline{\delta_{i,j}}$} \\ \hline
			\begin{tabular}[c]{@{}l@{}}number \\ of nodes\end{tabular} &
			\multicolumn{1}{l}{Linear $L_{\infty}$} &
			\multicolumn{1}{l|}{Quad $L_{\infty}$} &
			\multicolumn{1}{l}{{Linear $L_{\infty}$}} &
			\multicolumn{1}{l|}{{Quad $L_{\infty}$}} &
			\multicolumn{1}{l}{{ Linear $L_{\infty}$}} &
			\multicolumn{1}{l}{{ Quad $L_{\infty}$}} \\ \hline
			41$\times$41   & 2.047e-5      & 6.792e-7     & 7.741e-6    & 1.704e-5   & 1.957e-5     & 8.614e-7     \\
			81$\times$81   & 6.198e-6      & 1.722e-7     & 2.827e-6    & 4.256e-6   & 6.039e-6     & 1.652e-7     \\
			151$\times151$ & 1.7752e-6      & 4.9798e-8      & 8.4515e-7    & 1.2271e-6   & 1.7526e-6     & 4.8523e-8      \\ \hline
		\end{tabular}
	}
\end{table}

We extended the above tests to the 3-D example, with the results summarized in Table~\ref{tab:3D rhs}. Under the first condition, where accurate RHS values were used, the quadratic scheme significantly outperformed the linear scheme—achieving substantially lower numerical errors with only $26^3$ nodes, compared to the linear scheme using $101^3$ nodes. However, this advantage was diminished under the second condition, where inaccurate RHS values $\overline{b_{i,j,k}}$ were applied to near-boundary nodes. In this case, the linear scheme once again benefited from the altered RHS values, while the quadratic scheme exhibited a marked degradation in accuracy, more pronounced than what was observed in the 2-D case. Under the third condition, where the inaccurate RHS values were calibrated using the approximated coefficients $\overline{\delta_{i,j,k}}$ (i.e., applying $\widetilde{b_{i,j,k}} = \overline{b_{i,j,k}} / \overline{\delta_{i,j,k}}$), both the linear and quadratic schemes produced results closely matching those of the first condition. Notably, the quadratic scheme even yielded slightly improved accuracy compared to the original condition with exact RHS values, likely due to truncation error cancellation in this particular configuration.

These findings suggest that the use of approximated $\overline{\delta_{i,j}}$ (or $\overline{\delta_{i,j,k}}$ for 3-D problems) is a practical and effective approach for compensating the effects of RHS inaccuracies in the quadratic scheme, particularly in scenarios where accurate RHS data is not accessible near boundaries. Such situations frequently arise in plasma simulations conducted using the particle-in-cell (PIC) method.
\begin{table}[]
\caption{3-D test case with nonuniform $b_{i,j,k}$---the numerical error results obtained using different treatments of RHS values.}
\label{tab:3D rhs}
\setlength{\extrarowheight}{5pt} 
\resizebox{\textwidth}{!}{%
\begin{tabular}{l|cc|cc|cc}
\hline
        & \multicolumn{2}{c|}{Accurate RHS $b_{i,j,k}$} & 
        \multicolumn{2}{c|}{Inaccurate RHS $\overline{b_{i,j,k}}$} & 
        \multicolumn{2}{c}{Modified $\widetilde{b_{i,j,k}}=\overline{b_{i,j,k}}/\overline{\delta_{i,j,k}}$} \\ \hline
\begin{tabular}[c]{@{}l@{}}number \\ of nodes\end{tabular} &
  \multicolumn{1}{l}{Linear $L_{\infty}$} &
  \multicolumn{1}{l|}{Quad $L_{\infty}$} &
  \multicolumn{1}{l}{{Linear $L_{\infty}$}} &
  \multicolumn{1}{l|}{{Quad $L_{\infty}$}} &
  \multicolumn{1}{l}{{ Linear $L_{\infty}$}} &
  \multicolumn{1}{l}{{ Quad $L_{\infty}$}} \\ \hline
$26^3$   & 2.007e-4      & 1.113e-5     & 7.294e-5    & 1.768e-4   & 1.880e-4     & 6.039e-6     \\
$51^3$   & 4.939e-5      & 2.348e-6     & 2.076e-5    & 4.621e-5   & 4.779e-5     & 1.747e-6     \\
$101^3$ & 1.284e-5      & 5.281e-7      & 5.724e-6    & 1.187e-5   & 1.263e-5     & 4.599e-7      \\ \hline
\end{tabular}
}
\end{table}

\section{Conclusions}
In this paper, we investigated the numerical errors arising from inaccurate right-hand side (RHS) values near irregularly shaped Dirichlet boundaries when solving the electrostatic Poisson equation using the embedded finite difference method. Both linear and quadratic boundary treatments were considered. Such RHS inaccuracies commonly occur in plasma flow simulations using the particle-in-cell (PIC) method.

In the 1-D setting, we extended the analysis of Jomaa and Macaskill\cite{jomaa2005embedded} to derive explicit expressions for the numerical error as a function of local truncation errors at all nodes, incorporating the effects of RHS inaccuracies. This derivation was based on the assumption of a locally uniform RHS distribution. Our results showed that the numerical error can be decomposed into two boundary-induced components, $\xi_i^L$ and $\xi_i^R$, and an interior contribution, $\xi_i^{in}$. For the linear scheme, the RHS inaccuracies modified the truncation errors at near-boundary nodes: although these errors remain zeroth-order, their magnitude was reduced. This reduction in truncation error directly translated to a smaller boundary-induced error component. Since the interior error component remained unchanged, this implies that the RHS inaccuracies actually improved the solution quality of the linear scheme—an unexpected outcome. In contrast, for the quadratic scheme, RHS inaccuracies significantly worsened the truncation errors at near-boundary nodes by introducing an additional zeroth-order term, which in turn degraded the overall solution quality by increasing the boundary-induced error component. With this insight, the quadratic scheme’s error structure becomes more similar to that of the linear scheme, with boundary-induced errors dominating the global error. While our theoretical analysis assumed uniform RHS values, numerical experiments in 1-D confirmed these findings. Notably, as the grid is refined, the actual $\delta$ values—defined as the ratio between inaccurate and accurate RHS values—converge to the approximated $\overline{\delta}$ values derived under the uniform RHS assumption.

For the 2-D case, we evaluated how the magnitude of truncation errors at near-boundary nodes changes when inaccurate RHS values are used. To facilitate the analysis, we assumed the ``average" Dirichlet boundary shape represented by a straight line that's uniquely determined by $(\theta_x,\theta_y)$ pair associated to a near-boundary node, where $\theta_x$ and $\theta_y$ denotes the normalized distance to the boundary along axial directions. Our results showed that, for the linear scheme, the inaccurate RHS values generally reduced the truncation error near the boundary---except in extreme $(\theta_x, \theta_y)$ cases. For the quadratic scheme, however, the inaccurate RHS values significantly degraded the truncation errors, reducing them from first-order to zeroth-order accuracy---similar to the behavior observed in the lD setting. A comparison of the two schemes under RHS inaccuracies revealed that the linear scheme consistently outperformed the quadratic one, in the sense that the truncation errors at near-boundary nodes were smaller across nearly all $(\theta_x,\theta_y)$ combinations. Since a complete analytical expression for the numerical error was not available in 2D or 3D, we conducted numerical experiments in both settings. The results confirmed that modifications in truncation error near the boundary directly influence the global numerical error, and further supported that RHS inaccuracies tend to improve the linear scheme while degrading the performance of the quadratic scheme.

Given the severe degradation observed in the quadratic scheme under RHS inaccuracies, we also explored a mitigation strategy: replacing the inaccurate RHS values at near-boundary nodes with corrected values using the approximated ratios $\overline{\delta_{i,j}}$ (and $\overline{\delta_{i,j,k}}$ in 3D). Our numerical results demonstrated that this approach is both simple and highly effective in restoring solution accuracy for the quadratic scheme.

These findings offer valuable guidance for improving the accuracy and robustness of plasma simulations using the PIC method. We recommend adopting the linear embedded scheme as the default choice---not only because it preserves the symmetry structure favored by many efficient solvers, but also because it exhibits a counterintuitive benefit: the presence of RHS inaccuracies---commonplace in practical PIC simulations---can actually improve the solution accuracy. Moreover, numerical tests revealed that the linear scheme is more robust in a practical sense: even under increased grid resolution, it consistently maintained accuracy. In contrast, the quadratic scheme experienced sudden and severe accuracy degradation at specific grid sizes (e.g., $161 \times 161$ in our 2D tests), particularly when some nodes became trapped between two Dirichlet interfaces---an edge case to which the linear scheme was far less sensitive.If one still prefers to use the quadratic scheme due to its  accuracy advantages under ideal RHS conditions, the proposed $\overline{\delta}$-based correction strategy provides a practical and effective remedy for its sensitivity to RHS errors.

\section*{Acknowledgments}
We thank Professor Macaskill for his invaluable clarifications regarding aspects of his error analysis work. This work was supported by the National Natural Science Foundation of China (Grant No.12405245), Hunan Provincial Natural Science Foundation (Grant No.2024JJ6329), and Excellent Young Scholars Project of Hunan Provincial Department of Education (Grant No.23B0665).

\section*{Data Availability Statement}
All data generated or analyzed during this study are included within the manuscript. Numerical experiments can be reproduced using the described algorithms and parameters.








  \bibliographystyle{elsarticle-num} 
  \bibliography{JCP-Kai}






\end{document}